\documentclass[12pt,a4paper, twoside]{article}

\usepackage{amsmath,amssymb,amsthm,amsfonts,mathrsfs,amscd,environ}
\usepackage{dsfont}
\usepackage{latexsym,enumerate,geometry,extarrows}
\usepackage{stmaryrd}

\usepackage{amsfonts}
\usepackage{xcolor}

\usepackage{verbatim,fancyhdr,enumitem}%for comment paragraph fancyhdr,
 \newtheorem{thm}{Theorem}[section]
 \newtheorem{lem}[thm]{Lemma}
 
 \newtheorem{exa}[thm]{Example}
 
 \newtheorem{defn}{Definition}[section]
 \newtheorem{rem}{Remark}[section]
 
 \numberwithin{equation}{section}

\def\dif{{\mathord{{\rm d}}}}
\def\no{\nonumber}

\def\mR{{\mathbb R}}
\def\mE{{\mathbb E}}
\def\mM{{\mathbb M}}

\def\mE{{\mathbb E}}
\def\mF{{\mathbb F}}

\def\mJ{{\mathbb J}}

\def\mM{{\mathbb M}}

\def\mR{{\mathbb R}}
\def\mS{{\mathbb S}}

\def\sF{{\mathscr F}}

\def\sB{{\mathscr B}}

\def\sF{{\mathscr F}}

\def\sM{{\mathscr M}}

\def\sP{{\mathscr P}}

\def\bd{\begin{defn}}
\def\ed{\end{defn}}
\def\bl{\begin{lem}}
\def\el{\end{lem}}
\def\bt{\begin{thm}}
\def\et{\end{thm}}
\def\br{\begin{rem}}
\def\er{\end{rem}}
\allowdisplaybreaks
\title{{\bf Extended MF-FBSDEs with nonlinear domination-monotonicity  conditions and stochastic optimal controls  of  Linear System with  quadruple controls}%\footnote{*.}
}

\author{
{\bf Hao Wu$^{a)}$ }\\
\footnotesize{$^{a)}$School of Mathematics and Statistics,
South-Central Minzu University}\\
\footnotesize{ Wuhan, Hubei 430000, P.R.China}\\
\footnotesize{Email: wuhaomoonsky@163.com},
}

\begin{document}

\maketitle

\begin{abstract}
This paper extends the domination-monotonicity conditions, which guarantee the well-posedness of extended mean-filed forward-backward stochastic differential equations (extended MF-FBSDEs), from the previously studied linear framework to a nonlinear setting by incorporating nonlinear adjoint functions. Utilizing this generalized well-posedness result for extended MF-FBSDEs in conjunction with other refined analytical techniques, we address two classes of stochastic quadruple optimal controlled  problems: a linear-convex   problem and a linear-quadratic   problem with input constraints that are permitted to be time-dependent and random. For each problem, we establish the existence and uniqueness of optimal controls and derive their explicit closed-form representations.
\end{abstract}\noindent

\noindent

  Quadruple  controls; Linear-convex problem.

%\vskip 2cm
\section{Introduction}
Throughout this paper,  let  $(\Omega, \sF, \mF, P)$  be  a complete probability space with filtration $\mF= \{\sF_{t}\}_{t\geq 0}$ satisfying the usual conditions(i.e., it is increasing and right continuous,  $\sF_{0}$ contains all $P$-null sets) taking along
 a standard  $d$-Brownian motion $W(t):=(W_{1}(\cdot),W_{2}(\cdot), \cdots, W_{d}(\cdot) )^{\top},$ where the superscript $"\top"$ represents the transpose of a vector or a matrix.

The pervasive presence of uncertainty is a fundamental challenge in modeling and managing complex systems across various fields, from finance and engineering to biology and logistics. Stochastic control theory provides a powerful mathematical framework for addressing this challenge, enabling the design of strategies that are robust to random fluctuations. At its core, a stochastic control problem involves making sequential decisions over time to optimize a performance criterion-such as maximizing profit or minimizing cost when the underlying system dynamics are subject to random noise, often modeled by stochastic differential equations.

In stochastic control theory, the quest for optimal strategies governing dynamical systems subject to uncertainty finds a powerful mathematical representation in forward-backward stochastic differential equations(FBSDEs), which focuses on the maximum principle, variational techniques and a dual perspective are applied to derive necessary conditions for open-loop optimal control. This leads to a coupled FBSDE, often referred to in control theory as a stochastic Hamiltonian system (cf. \cite{BKMM,LX2,TianYu2023,WangZhang2017,XuXieZhang2017,Yong2010} ).  These systems inherently split the problem: a forward SDE dictates the evolution of the state, and a backward SDE, linked by a terminal condition, propagates the cost or adjoint variables. This formulation is crucial, as it directly leads to the stochastic maximum principle-a cornerstone for establishing optimality. Thus, FBSDEs are not merely a descriptive tool but a fundamental construct for solving, analyzing, and approximating solutions to a wide array of stochastic control problems (cf.\cite{AnderssonAnderssonOosterlee2023,Kushner1990, Yu2022}).  In the subsequent phase, issues such as existence and uniqueness of the FBSDE are examined(cf.\cite{Antonelli1993,CarmonaDelarue2013,HuPeng1995,MaProtterYong1994}), and from these findings, corresponding properties like existence and uniqueness for optimal control are established (as seen in works like \cite{BKMM1,ExarchosTheodorou2018,LiuNiuWangYu2026}).

Basing on the theory of FBSDEs,  particular significance produces  the theory of stochastic control produces the framework of mean field forward-backward stochastic differential equations (MF-FBSDEs). The existence and uniqueness results are investigated with different conditions(cf. \cite{BayraktarZhang2023,BensoussanYamZhang2015,BuckdahnDjehicheLiPeng2009,BuckdahnLiPeng2009,CarmonaDelarue2013,CarmonaDelarue2015}). This framework provides a powerful and elegant mathematical tool for analyzing stochastic control problems involving a large population of interacting agents, where the dynamics of each individual are influenced not only by its own state and control but also by the statistical distribution of the entire population's states and controls. The forward equation describes the evolution of the state process, while the backward equation inherently encodes the necessary optimality conditions derived from the stochastic maximum principle. The "mean-field" component elegantly captures the aggregate effect of the population, allowing the complex, high-dimensional problem to be approximated by a simpler, centralized control problem for a representative agent. Consequently,  MF-FBSDEs have found profound applications in areas such as mean-field games, systemic risk modeling in large financial networks, and optimal resource allocation in massive systems, offering a tractable pathway to derive decentralized strategies and analyze their limiting behaviors(cf.\cite{Ahuja2016,Ahuja2019,Alasseur2020,BayraktarZhang2023,DFF1,DFF2,Djete2023,LasryLions2007} ).

Especially, studying linear-quadratic (LQ) optimal control problems for such mean-field FBSDEs is of particular theoretical and practical significance(cf.\cite{LiSunXiong2019,TianYu2023}).
The motivation for investigating LQ problems in this context is twofold. Theoretically , the LQ framework offers a tractable yet rich structure that often admits explicit or computationally accessible solutions. It serves as a crucial benchmark and a foundational building block for understanding more general nonlinear problems. By examining mean-field FBSDEs through the LQ problems, one can derive sharp conditions for solvability, obtain closed-form optimal controls, and gain profound insight into the interplay between individual optimization, stochasticity, and mean-field coupling. Practically , LQ mean-field control models are directly applicable to numerous domains. Examples include portfolio optimization in finance with price impacts, consensus control in multi-agent robotic systems, and macroeconomic models where agents' decisions depend on aggregate economic indicators(cf. \cite{Bardi2012,LiLiYu2020,LinJiangZhang2019}). Therefore, the study of linear-quadratic optimal control for mean-field FBSDEs not only advances the mathematical theory of stochastic control and mean-field games but also provides a versatile toolkit for designing and analyzing optimal strategies in complex, interconnected systems.

In this article,  we intend to investigate the following    extended MF-FBSDEs with nonlinear domination-monotonicity  conditions by  setting\\ $\theta_{i}(t):=((\mE[X_{1}(t)])^{\top}, (\mE[Y_{1}(t)])^{\top}, (\mE[X_{2}(t)])^{\top}, (\mE[Y_{2}(t)])^{\top},   X_{i}(t)^{\top}, Y_{i}(t)^{\top}, Z_{i}(t)^{\top})^{\top}, i=1,2 : $
\begin{align}\label{1}
\begin{cases}
&\dif X_{1}(t)=b_{1}(t,\theta_{1}(t))\dif t
+\sum^{d}_{j=1}\sigma_{1j}(t,\theta_{1}(t))\dif W_{j}(t), t\in [0,T]\\
&\dif Y_{1}(t) =f_{1}(t,\theta_{1}(t))\dif t+\sum^{d}_{j=1} Z_{1j}(t)\dif W_{j}(t)\\
& X_{1}(0)=\Psi_{1}(Y_{1}(0),Y_{2}(0)),  Y_{1}(T)=\Phi_{1}(X_{1}(T),X_{2}(T) ),
\end{cases}
\end{align}

 \begin{align}\label{2}
\begin{cases}
&\dif X_{2}(t)=b_{2}(t,\theta_{2}(t))\dif t
 +\sum^{d}_{j=1}\sigma_{2j}(t,\theta_{2}(t))\dif W_{j}(t), t\in [0,T]\\
&\dif Y_{2 }(t) =f_{2}(t,\theta_{2}(t))\dif t+\sum^{d}_{j=1} Z_{2j}(t)\dif W_{j}(t)\\
& X_{2 }(0)=\Psi_{2}(Y_{1}(0),Y_{2}(0)),  Y_{2}(T)=\Phi_{2}(X_{1}(T),X_{2}(T) ).
\end{cases}
\end{align}
Denote $V_{i}(\cdot):=(X_{i}(t)^{\top}(\cdot),$ $Y_{i}(\cdot)^{\top},Z_{i}(\cdot)^{\top})^{\top}, i=1,2,$   which are the unknown processes   with    $Z_{i}(\cdot):=(Z_{i1}(\cdot)^{\top},Z_{i2}(\cdot)^{\top},\cdots,Z_{id}(\cdot)^{\top})^{\top}, \sigma_{i}(\cdot):=(\sigma_{i1}(\cdot)^{\top},\sigma_{i2}(\cdot)^{\top},\cdots,\sigma_{id}(\cdot)^{\top})^{\top},  i=1,2,.$   The coefficients $\Phi_{i}: \mR^{n}\times \mR^{n}\rightarrow \mR^{n}, $ $\Psi_{i}: \mR^{n}\times \mR^{n}
  \rightarrow \mR^{n} $  and $f_{i},b_{i}, \sigma_{ij}: \Omega \times [0,T] \times \mR^{n+n+n+n+n+n+nd}\rightarrow \mR^{n}, j=1,2,\cdots d, i=1,2.$ For simplicity, denote $\Gamma_{i}(\cdot):=(f_{i}(\cdot)^{\top},b_{i}(\cdot)^{\top},\sigma_{i}(\cdot)^{\top})^{\top}, i=1, 2. $   Besides,  the combination of Eq.\eqref{1} and Eq.\eqref{2} is denoted by $(\pi).$  Consequently, all coefficients of the System $(\pi)$ is encompassed by $(\Phi_{i}, \Psi_{i}, \Gamma_{i}), i=1,2$. Obviously, there exist two initial couplings in the System $(\pi).$    Later,  when we  study the optimal control problems associated with System $(\pi),$  two initial values will appear served as the two initial controls. This article will establish the results regarding the well-posedness of System $(\pi),$ specifically addressing existence, uniqueness, and  the applications in stochastic optimal control problems.

 The motivation for this work stems from two aspects.  The first motivation follows from the work of \cite{BuckdahnLiPeng2009}. In   \cite{BuckdahnLiPeng2009},  Buckdahn et al. investigate the following  mean field FBSDE:
  \begin{align}\label{3}
\begin{cases}
&\dif X_{1}(t)=\mE[b(t,X_{1}(t),Y_{1}(t),\mathbf{u})]|_{\mathbf{u}=(X_{1}(t),Y_{1}(t), Z_{1}(t) )}\dif t\\
&\quad\quad\quad\quad +\sum^{d}_{j=1}\mE[\sigma_{j}(t,X_{1}(t),Y_{1}(t),\mathbf{u})]|_{\mathbf{u}=(X_{1}(t),Y_{1}(t), Z_{1}(t) )}\dif W_{j}(t), t\in [0,T]\\
&\dif Y_{1}(t) =\mE[f(t,X_{1}(t),Y_{1}(t), \mathbf{u})]|_{\mathbf{u}=(X_{1}(t),Y_{1}(t),Z_{1}(t)  )}\dif t +\sum^{d}_{j=1} Z_{1j}(t)\dif W_{j}(t)\\
& X_{1}(0)=x,  Y_{1}(T)=\Phi_{1}(X_{1}(T),X_{2}(T) ).
\end{cases}
\end{align}
 and
 \begin{align}\label{3+}
\begin{cases}
&\dif X_{2}(t)=\mE[b(t,X_{1}(t),Y_{1}(t),\mathbf{u})]|_{\mathbf{u}=(X_{2}(t),Y_{2}(t), Z_{2}(t) )}\dif t\\
&\quad\quad\quad\quad +\sum^{d}_{j=1}\mE[\sigma_{j}(t,X_{1}(t),Y_{1}(t),\mathbf{u})]|_{\mathbf{u}=(X_{2}(t),Y_{2}(t), Z_{2}(t) )}\dif W_{j}(t), t\in [0,T]\\
&\dif Y_{2}(t) =\mE[f(t,X_{1}(t),Y_{1}(t), \mathbf{u})]|_{\mathbf{u}=(X_{2}(t),Y_{2}(t),Z_{2}(t)  )}\dif t +\sum^{d}_{j=1} Z_{2j}(t)\dif W_{j}(t)\\
& X_{2}(0)=x,  Y_{2}(T)=\Phi_{2}(X_{1}(T),X_{2}(T) ).
\end{cases}
\end{align}
Under the certain  conditions,  they got the existence and uniqueness result of solution. Moreover, the authors prove that the unique solution for the above system  is closely related to  a nonlocal partial differential equation. Compare with the Eqs.\eqref{3}, \eqref{3+},  we will give several different points. The first  is that   the initial values and the terminal values of  the system $(\pi)$ contain two couplings. As far as we known, previous MF-FBSDEs did not take into account this double coupled  duality.  Later,  when we  study the optimal control problem associated with system $(\pi),$  two initial values will appear served as the two initial controls,  which is one of our contributions.  The second point is that the coefficients satisfy nonlinear
domination-monotonicity conditions but not linear domination-monotonicity. The last point is that we consider an extended  coupling compared with the case in \eqref{3}, \eqref{3+}.

The other  motivation  is that  Tian and Yu \cite{TianYu2023}  and    Liu et al. \cite{LiuNiuWangYu2026}  studied  the single initial value control problems. Compared with the systems in \cite{TianYu2023} and \cite{LiuNiuWangYu2026}.  We give several different points.  First, we consider double initial value control problems. In subsequent control applications for such problems, there will be two evolution paths, where the first path influences the evolution of the second, leading to an essential difference in the value function.  Second, compared with the system in Tian and Yu \cite{TianYu2023},  we consider a more generalized mean field system. Furthermore,  there exist two initial couplings while Liu et al. \cite{TianYu2023} didn't consider initial couplings and  Liu et al. \cite{LiuNiuWangYu2026} only consider one initial couplings. Third, inspired by  \cite{LiuNiuWangYu2026}, we consider the case of nonlinear
domination-monotonicity conditions  on coefficients while Tian and Yu \cite{TianYu2023} didn't consider the  case of nonlinear
domination-monotonicity conditions for the mean field system.     In detail,  the controlled system is modeled by the following linear coupled SDE defined on $[0, T],$
\begin{align}\label{4}
\begin{cases}
&\dif X_{1}(t)=[A_{1}(t)X_{1}(t) + \bar{A}_{2}(t)\mE[X_{2}(t)] +B_{1}(t)u_{1}(t)+\tau\bar{B}_{1}(t)\mE[u_{1}(t)]\\
&\quad\quad\quad\quad\quad\quad\quad\quad\quad\quad\quad\quad\quad\quad\quad\quad+\bar{A}_{1}(t)\mE[ X_{1}(t)]+\rho_{1}(t)]\dif t\\
&\quad\quad\quad\quad+\sum^{d}_{j=1}[C_{1j}(t)X_{1}(t) +D_{1j}(t)u_{1}(t)+\kappa_{1j}]\dif W_{j}(t), t\in [0,T]\\
&X_{2}(t)=[A_{2}(t)X_{2}(t)+B_{2}(t)u_{2}(t)+\tau\bar{B}_{2}(t)\mE[u_{2}(t)]+\bar{A}_{2}(t)\mE[ X_{1}(t)] +\rho_{2}(t) ]\dif t\\
&\quad\quad\quad\quad +\sum^{d}_{j=1}[C_{2j}(t)X_{2}(t) +D_{2j}(t)u_{2}(t)+\kappa_{2j}]\dif W_{j}(t), t\in [0,T],\\
& X_{1}(0)=H\xi_{1}+x_{0}, X_{2}(0)=H\xi_{2}+x_{0}.
\end{cases}
\end{align}
where $\tau>0$ is a sufficiently small constant, $\xi_{i}\in \mR^{m}, \rho_{i}(\cdot)\in L^{2}_{\mF}(\mR^{n}), \kappa_{ij}\in L^{2}_{\mF}(\mR^{n}) ,i=1,2, j=1,2,\cdots, d,$       $H$ is an   suitable matrix and $A_{i}(\cdot), B_{i}(\cdot),\bar{B}_{i}(\cdot), \bar{A}_{i}(\cdot), C_{ij}(\cdot), D_{ij}(\cdot), i=1,2, j=1,2,\cdots, d,$ are suitable matrix-valued function. The detailed definitions  will be provided below.  In contrast to the most of the literature where the control system only involves the process control $u_{1}(\cdot),$  The System \eqref{4} not only involves two the process control $u_{1}(\cdot), u_{2}(\cdot),$  but also involves double initial controls $\xi_{1},\xi_{2}$, which is called stochastic quadruple optimal control problems. Compared with the single initial control in \cite{LiuNiuWangYu2026}, we consider the more generalize  double initial controls, which maybe have more practical significance.  This system leads to the associated Hamiltonian  system having more complicated coupling.

$\mathbf{Linear-convex\, problem}$($\mathbf{Problem\,(LC)}$):  when we study this linear-convex problem,  compared to conventional convex control problems,  our value functional  has the following form($\tau >0$ is a sufficiently small constant):
\begin{align}\label{5}
J(\xi_{1},\xi_{2}, u_{1}(\cdot),u_{2}(\cdot))&:=f_{11}(\xi_{1}+\tau\xi_{2})+ f_{12}(\xi_{2}+\tau\xi_{1})\no\\
&+\mE\bigg\{\bigg.f_{21}(X_{1}(T)+X_{2}(T) ) +f_{22}( X_{1}(T)+X_{2}(T) ) \no \\
 &  + \int^{T}_{0}f_{31}(t,X_{1}(t) )\dif t +\int^{T}_{0}f_{32}(t,X_{2}(t))\dif t\no\\
 &+\int^{T}_{0}f_{41}(t,u_{1} (t) )\dif t +\int^{T}_{0}f_{42}(t,  u_{2}(t))\dif t \bigg\}\bigg..
\end{align}
where $f_{2i}: \Omega \times \mathbb{R}^n \rightarrow [0,\infty),$  $f_{3i}: \Omega \times [0,T] \times \mathbb{R}^n \rightarrow [0,\infty)$  and $f_{1i}: \mathbb{R}^m \rightarrow [0, \infty)$ and $f_{4i}: \Omega \times [0, T] \times \mathbb{R}^k \rightarrow [0, \infty).$  We  require $f_{ji}, j=1,2,3,4, i=1,2$ to satisfy  convexity,  which will be presented  later.  We intend to find a quartet of $(\eta^{*},\xi^{*},u_{1}^{*}(\cdot),u_{2}^{*}(\cdot))$ such that
\begin{align}\label{5+}
J(\xi^{*}_{1},\xi^{*}_{2},u_{1}^{*}(\cdot),u_{2}^{*}(\cdot)):=\inf_{(\xi_{1},\xi_{2}, u_{1}(\cdot),u_{2}(\cdot))\in \mR^{m}\times \mR^{m} \times L^{2}_{\mF}(\mR^{k})\times L^{2}_{\mF}(\mR^{k})}J(\xi_{1},\xi_{2}, u_{1}(\cdot),u_{2}(\cdot)).
\end{align}
We establish the existence and uniqueness of optimal controls and derive their explicit closed-form representations. Further details for linear-convex problem will be elaborated in the following sections.

$\mathbf{Linear-quadratic\, problem\, with\, input\, constraints}\,( \mathbf{Problem\,(LQ-IC)}):$ In this case, we retain the linear controlled system \eqref{4} and reformulate the performance criterion \eqref{5} into the following quadratic form:

\begin{align}\label{6}
 &\mJ(\xi_{1},\xi_{2}, u_{1}(\cdot),u_{2}(\cdot)):=\frac{1}{2}\langle M_{1}\xi_{1}, \xi_{1}\rangle+\frac{1}{2}\langle M_{2}\xi_{2}, \xi_{2}\rangle,\no\\
 &+\frac{1}{2}\mE\bigg\{\bigg. \langle G_{1}(X_{1}(T)+X_{2}(T)), X_{1}(T)+X_{2}(T)\rangle + \langle G_{2}(X_{1}(T)+X_{2}(T)), X_{1}(T)+X_{2}(T)\rangle  \no \\
 &  + \int^{T}_{0}\langle Q_{1}(t)X_{1}(t), X_{1}(t)\rangle\dif t+ \int^{T}_{0}\langle Q_{2}(t)X_{2}(t), X_{2}(t)\rangle\dif t\no\\
 &+\int^{T}_{0}\langle R_{1}(t) u_{1}(t), u_{1}(t)\rangle\dif t +\int^{T}_{0}\langle R_{2}(t)u_{2}(t), u_{2}(t)\rangle\dif t \bigg\}\bigg.,
\end{align}
where $M_{1}(\cdot),M_{2}(\cdot),G _{1}(\cdot),G _{2}(\cdot), Q_{1}(\cdot) ,Q_{2}(\cdot), R_{1}(\cdot), R_{2}(\cdot)$ are suitable matrix-valued functions, which will be presented in Section 5.
Unlike Problem (LC), we now constrain the control variables (inputs) $\xi_{i},i=1,2$ and $u_{i}(\cdot),i=1,2$ to reside within nonempty, closed, convex constraint sets:
\begin{equation}\label{as1}
    U_0 \subset \mathbb{R}^m \text{ and } U(\cdot) \equiv \{U(\omega,t) \subset \mathbb{R}^k, (\omega,t) \in \Omega \times [0, T]\}
   \end{equation}
satisfying certain conditions, respectively. Let
\begin{equation}\label{as}
    \mathcal{U} := \{u(\cdot) \in L^2_{\mF}(\mathbb{R}^k) \mid u(\omega,t) \in U(\omega, t) \text{ for almost all } (\omega, t) \in \Omega \times [0,T]\}
    \end{equation}
and refer to $U_0\times U_0 \times \mathcal{U}\times \mathcal{U}$ as the admissible control set. Our linear-quadratic stochastic optimal control problem with input constraints (IC) is presented as follows.

\textbf{Problem (LQ-IC):} Determine a quartet of admissible control inputs $(\xi^*_{1},\xi^*_{2}, u^*_{1}(\cdot),u^*_{2}(\cdot))$ such that
\begin{equation}
    \mJ(\xi^*_{1},\xi^*_{2}, u^*_{1}(\cdot),u^*_{2}(\cdot))= \inf_{(\xi_{1},\xi_{2}, u_{1}(\cdot),u_{2}(\cdot)) \in U_0\times U_0 \times \mathcal{U}\times \mathcal{U}} \mJ(\xi_{1},\xi_{2}, u_{1}(\cdot),u_{2}(\cdot)).
    \label{eq:9}
\end{equation}

In this scenario, $(\xi^*_{1},\xi^*_{2}, u^*_{1}(\cdot)), u^*_{2}(\cdot))$, $X^*_{i}(\cdot) := x(\cdot; \xi^*_{i}, u^*_{i}(\cdot)), i=,2,$, and $(\xi^*_{1},\xi^*_{2}, u^*_{1}(\cdot),$ $u^*_{2}(\cdot), X^{*}_{1}(\cdot), X^{*}_{2}(\cdot) )$ are designated as a  quartet of optimal controls, the associated optimal state, and an optimal sextet for Problem (LQ-IC), respectively.

\section{Notations}
  For $x, y \in \mR^{n},$ we use $|x |$  to denote the Euclidean norm of
$x,$  $\langle x, y\rangle$ to denote the Euclidean inner product.
 For $A\in \mR^{ n\times d},$     $|A |$ represents  $\sqrt{\mathrm{Tr} (AA^{\top})}.$    Let $\sP$ stand for the $\mathbb{F}$-progressively measurable $\sigma$-field and
$\sB(\mathbb{R}^n)$ stand for the Borel $\sigma$-field on $\mathbb{R}^n$. $\mS^{n}\subset \mR^{n\times n}$ represents the set of all the symmetrical matrices.
Next, We intend to define the following Banach spaces of random vectors or stochastic processes $(p\geq 2).$
\begin{itemize}
\item[1)] $L^{2}_{\sF_{T}}(\mR^{n})$ is the set of $\sF_{T}-$measuurable random vectors $\zeta: \Omega\rightarrow \mR^{n}$ such that $\|\zeta\|_{L^{2}_{\sF_{T}}(\mR^{n})}:=\{\mE[|\zeta|^{2}]\}^{\frac{1}{2}}<\infty.$
\item[2)] $L^{\infty}_{\sF_{T}}(\mR^{n})$ is the set of $\sF_{T}-$measuurable random vectors $\zeta: \Omega\rightarrow \mR^{n}$ such that $\|\zeta\|_{L^{\infty}_{\sF_{T}}(\mR^{n})}:=\operatorname{ess sup}_{\omega\in \Omega}|\zeta|<\infty.$
\item[3)] $L^{2}_{\mF}(\mR^{n})$ is the set of $\sP-$measurable stochastic processes $\varphi: \Omega\times [0,T]\rightarrow \mR^{n}$ such that $\|\varphi\|_{L^{2}_{\mF}(\mR^{n})}:=\bigg\{\bigg. \mE\int^{T}_{0}|\varphi(t)|^{2}\dif t   \bigg\}\bigg.^{\frac{1}{2}}<\infty.$

\item[4)] $S^{2}_{\mF}(\mR^{n})$ is the set of   continuous  stochastic processes $\varphi: \Omega\times [0,T]\rightarrow \mR^{n}$ such that $\|\varphi\|_{S^{2}_{\mF}(\mR^{n})}:=\bigg\{\bigg. \mE[\sup_{0\leq t\leq T}|\varphi(t)|^{2}] \bigg\}\bigg.^{\frac{1}{2}}<\infty.$

\item[5)] $L^{\infty}(\mR^{n})$  consists of all deterministic function such that $\|\psi\|_{L^{\infty}(\mR^{n})}:=\sup_{ t\in [0,T]}|\psi|<\infty.$

\end{itemize}

Furthermore, for the  sake of simplicity, we set
\begin{itemize}
\item[1)] $\sM_{\mF}(\mR^{n+n+nd}):=L^{2}_{\mF}(\mR^{n})\times L^{2}_{\mF}(\mR^{n})\times L^{2}_{\mF}(\mR^{nd}) $ equipped with the norm $\|\alpha(\cdot)\|_{\sM(\mR^{n+n+nd})}:=\{\mE[\Xi_{\alpha}]\}^{\frac{1}{2}}$ for any $\alpha(\cdot):=(\varphi_{1}(\cdot)^{\top}, \varphi_{2}(\cdot)^{\top},\varphi_{3}(\cdot)^{\top})^{\top}\in \sM_{\mF}(\mR^{n+n+nd})$ where
\begin{align}\label{99}
\Xi_{\alpha}:=\int^{T}_{0}|\varphi_{1}(t)|^{2}\dif t  +\int^{T}_{0}|\varphi_{2}(t)|^{2}\dif t+\int^{T}_{0}|\varphi_{3}(t)|^{2}\dif t.
\end{align}

\item[2)]$\mM_{\mF}(\mR^{n+n+nd}):=S^{2}_{\mF}(\mR^{n})\times S^{2}_{\mF}(\mR^{n})\times L^{2}_{\mF}(\mR^{nd}) $ equipped with the norm $\|V(\cdot)\|_{\mM(\mR^{n+n+nd})}:=\{\mE[\Lambda_{V}]\}^{\frac{1}{2}}$ for any $V(\cdot):=(x(\cdot)^{\top}, y(\cdot)^{\top}, z(\cdot)^{\top})^{\top}\in \mM_{\mF}(\mR^{n+n+nd})$ where
\begin{align}\label{301}
\Lambda_{V}:=\sup_{0\leq t\leq T}|x(t)|^{2}  +\sup_{0\leq t\leq T}|y(t)|^{2}\dif t+\int^{T}_{0}|z(t)|^{2}\dif t.
 \end{align}

\end{itemize}

\section{Well-posedness of  extended MF-FBSDE}
In this section, we examine the well-posedness of  System $(\pi)$.
 \subsection{Assumptions} To begin, we introduce the following assumptions on the coefficients.  The inspiration for the following  assumptions partially  comes from \cite{LiuNiuWangYu2026}.

$\mathbf{Assumption\, 1}:$
\begin{itemize}
\item[ (i)]$\Psi_{i},$ $ \Phi_{i}$ and $\Gamma_{i},i=1,2$  are each measurable relative to the $\sigma-$algebras $\sB(\mR^{n})\times \sB(\mR^{n}),$ $ \sF_{T}\times\sB(\mR^{n})\times \sB(\mR^{n})$ and $\sP\times\sB(\mR^{n+n+n+n+n+n+nd}),$ respectively. Furthermore, $\Phi_{i}(0,0)\in L^{2}_{\sF_{T}}(\mR^{n}), \Gamma_{i}(\cdot,0)\in \sM_{\mF}(\mR^{n+n+nd}),i=1,2,$  and there exists a positive constant $L$ such that for any $t\in [0,T], \omega\in \Omega, x', x''\in \mR^{n}, y', y''\in \mR^{n},y_{1}, y_{2}\in \mR^{n},$
   $$|\Psi_{1}(0,y_{2})|+|\Psi_{2}(y_{1},0)|+|l(t,x',y',x'',y'',0,0,0)|\leq L,$$
  where $l:=b_{i},\sigma_{i},f_{i},i=1,2.$

\item[(ii)](Lipschitz conditions) there exist two  constants $L_{b}>0,L_{\sigma}>0,L_{f}>0,L_{\Phi}>0 ,L_{\Psi}>0$ and two small enough constants $\varepsilon>0, \epsilon>0$ such that, for
$t\in [0,T],$ $x^{\prime},x^{\prime\prime}, \bar{x}^{\prime},\bar{x}^{\prime\prime}  ,x,\bar{x}, y^{\prime},$ $y^{\prime\prime}, \bar{y}^{\prime},\bar{y}^{\prime\prime},y,\bar{y}\in \mR^{n},$ \\ $z:=(z_{1}^{\top},z_{2}^{\top},\cdots,z_{d}^{\top})^{\top}\in \mR^{nd},\bar{z}:=(\bar{z}_{1}^{\top},\bar{z}_{2}^{\top},\cdots,\bar{z}_{d}^{\top})^{\top}\in \mR^{nd},  $
\begin{align*}
&|b_{i}(t,x^{\prime},y^{\prime},x^{\prime\prime},y^{\prime\prime},x,y,z )-b_{i}(t, \bar{x}^{\prime},\bar{y}^{\prime},\bar{x}^{\prime\prime},\bar{y}^{\prime\prime},\bar{x},\bar{y},\bar{z})|\\
&<\varepsilon|x^{\prime}-\bar{x}^{\prime}|+\varepsilon|y^{\prime}-\bar{y}^{\prime}|+\varepsilon|x^{\prime\prime}-\bar{x}^{\prime\prime}|
+\varepsilon|y^{\prime\prime}-\bar{y}^{\prime\prime}|+L_{b}|x-\bar{x}|+L_{b}|y-\bar{y}|+L_{b}|z-\bar{z}|,
\end{align*}

\begin{align*}
&|\sigma_{i}(t,x^{\prime},y^{\prime},x^{\prime\prime},y^{\prime\prime},x,y,z )-\sigma_{i}(t, \bar{x}^{\prime},\bar{y}^{\prime},\bar{x}^{\prime\prime},\bar{y}^{\prime\prime},\bar{x},\bar{y},\bar{z})|\\
&<\varepsilon|x^{\prime}-\bar{x}^{\prime}|+\varepsilon|y^{\prime}-\bar{y}^{\prime}|+\varepsilon|x^{\prime\prime}-\bar{x}^{\prime\prime}|
+\varepsilon|y^{\prime\prime}-\bar{y}^{\prime\prime}|+L_{\sigma}|x-\bar{x}|+L_{\sigma}|y-\bar{y}|+L_{\sigma}|z-\bar{z}|,
\end{align*}

\begin{align*}
&|f_{i}(t,x^{\prime},y^{\prime},x^{\prime\prime},y^{\prime\prime},x,y,z )-f_{i}(t, \bar{x}^{\prime},\bar{y}^{\prime},\bar{x}^{\prime\prime},\bar{y}^{\prime\prime},\bar{x},\bar{y},\bar{z})|\\
&<\varepsilon|x^{\prime}-\bar{x}^{\prime}|+\varepsilon|y^{\prime}-\bar{y}^{\prime}|+\varepsilon|x^{\prime\prime}-\bar{x}^{\prime\prime}|
+\varepsilon|y^{\prime\prime}-\bar{y}^{\prime\prime}|+L_{f}|x-\bar{x}|+L_{f}|y-\bar{y}|+L_{f}|z-\bar{z}|.
\end{align*}

\begin{align*}
|\Phi_{1}(x',x )-\Phi_{1}(\bar{x}',\bar{x})|\leq L_{\Phi}|x'-\bar{x}'|+\epsilon|x-\bar{x}|,\\
|\Phi_{2}(x',x )-\Phi_{2}(\bar{x}',\bar{x})|\leq\epsilon|x'-\bar{x}'|+L_{\Phi}|x-\bar{x}|.
\end{align*}

\begin{align*}
|\Psi_{1}(y',y )-\Psi_{1}(\bar{y}',\bar{y})|\leq L_{\Psi}|y'-\bar{y}'|+\epsilon|y-\bar{y}|,\\
|\Psi_{2}(y',y )-\Psi_{2}(\bar{y}',\bar{y})|\leq\epsilon|y'-\bar{y}'|+L_{\Psi}|y-\bar{y}|.
\end{align*}

\item[(iii)] There exist two constants $L_{1}>0,$ $L_{2}>0,$  a matrix $H\in \mR^{n\times m},$   several matrix-valued functions $B_{i}(\cdot),B_{i}(\cdot)\in L^{\infty}(\mR^{n\times k}), i=1,2$ and $D_{i}(\cdot)\in L^{\infty}(\mR^{n\times k})$ with $D_{i}(\cdot):=(D_{i1}(t)^{\top},D_{i2}(t)^{\top},\cdots,D_{di}(t)^{\top})^{\top},$  four $\sB(\mR^{m})-$measurable mappings $\bar{h}_{i1},\bar{h}_{i2}: \mR^{m} \rightarrow\mR^{m}, i=1,2$ and two $\sP\times \sB(\mR^{k})-$measurable mappings
    $h_{i}:\Omega\times [0,T]\times \mR^{k}\rightarrow \mR^{k}, i=1,2$ such that the following holds for $x^{\prime}_{i},x^{\prime\prime}_{i},y^{\prime}_{i},y^{\prime\prime}_{i},x_{i},y_{i},y\in \mR^{n},z_{i}\in \mR^{nd},$ $\bar{x}^{\prime}_{i},\bar{x}^{\prime\prime}_{i},\bar{y}^{\prime}_{i},\bar{y}^{\prime\prime}_{i},\bar{x}_{i},\bar{y}_{i}\in \mR^{n},\bar{z}_{i}\in \mR^{nd},i=1,2$.
\begin{itemize}
\item[1)]Adjoint function
\begin{align}\label{8}
\begin{cases}
&h_{i}(\cdot,0)\in L^{2}_{\mF}(\mR^{k}),i=1,2,\\
&|\bar{h}_{ik}( v_{1})-\bar{h}_{ik}( v_{2})|<L_{2}|v_{1}-v_{2}|,i=1,2,k=1,2,\\
&|h_{i}(t, u_{1})-h_{i}(t, u_{2})|<L_{2}|u_{1}-u_{2}|,i=1,2,\\
&\langle \bar{h}_{ik}(  v_{1})-\bar{h}_{ik}( v_{2}), v_{1}-v_{2} \rangle\leq -L_{3}|\bar{h}_{ik}( v_{1})-\bar{h}_{ik}( v_{2})|^{2},i=1,2,k=1,2,\\
&\langle h_{i}(t, u_{1})-h_{i}(t, u_{2}), u_{1}-u_{2} \rangle\leq -L_{3}|h_{i}( t,u_{1})-h_{i}( t,u_{2})|^{2},
\end{cases}
\end{align}
for any $(\omega,t)\in [0,T], v_{1},v_{2}, v'_{1},v'_{2},v\in \mR^{m},u_{1},u_{2}\in \mR^{k}.$

\item[2)]Domination conditions($i=1,2,$ $\tau$ is a sufficiently small positive constant):
\begin{align}\label{9}
\begin{cases}
&|\Psi_{1}(y_{1},y_{2})-\Psi_{1}(\bar{y}_{1},\bar{y}_{2})|\leq L_{2}|\bar{h}_{11}(\frac{H^{\top}y_{1}-\tau H^{\top}y_{2}}{1-\tau^{2}})-\bar{h}_{11}(\frac{H^{\top}\bar{y}_{1}-\tau H^{\top}\bar{y}_{2}}{1-\tau^{2}})|^{2},\\
&\quad \quad\quad\quad \quad\quad\quad \quad\quad\quad\quad+L_{2}|\bar{h}_{12}(\frac{H^{\top}y_{2}-\tau H^{\top}y_{1}}{1-\tau^{2}})-\bar{h}_{12}(\frac{H^{\top}\bar{y}_{2}-\tau H^{\top}\bar{y}_{1}}{1-\tau^{2}})|^{2},\\
&|\Psi_{2}(y_{1},y_{2})-\Psi_{2}(\bar{y}_{1},\bar{y}_{2})|\leq L_{2}|\bar{h}_{21}(\frac{H^{\top}y_{2}-\tau H^{\top}y_{1}}{1-\tau^{2}})-\bar{h}_{11}(\frac{H^{\top}\bar{y}_{2}-\tau H^{\top}\bar{y}_{1}}{1-\tau^{2}})|^{2},\\
&\quad \quad\quad\quad \quad\quad\quad \quad\quad\quad\quad+L_{2}|\bar{h}_{22}(\frac{H^{\top}y_{1}-\tau H^{\top}y_{2}}{1-\tau^{2}})-\bar{h}_{22}(\frac{H^{\top}\bar{y}_{1}-\tau H^{\top}\bar{y}_{2}}{1-\tau^{2}})|^{2},\\
&|l_{1}(t,x^{\prime}_{1},y^{\prime}_{1},x^{\prime\prime}_{1},y^{\prime\prime}_{1},x_{1},y_{1},z_{1})
-l_{1}(t,x^{\prime}_{1},y^{\prime}_{2},x^{\prime\prime}_{1},y^{\prime\prime}_{1},x_{1},y_{2},z_{2})|,\\
&\leq L_{2}(|h_{1}(t,B_{1}(t)^{\top}y_{1}+\tau \bar{B}_{1}(t)^{\top}y^{\prime}_{1}+D_{1}(t)^{\top}z_{1})\\
&\quad \quad\quad\quad \quad\quad\quad \quad\quad-h_{1}(t,B_{1}(t)^{\top}y_{2}+\tau \bar{B}_{1}^{\top}(t)^{\top}y^{\prime}_{2}+D_{1}(t)^{\top}z_{2}),\\
&|l_{2}(t,x^{\prime}_{1},y^{\prime}_{1},x^{\prime\prime}_{1},y^{\prime\prime}_{1},x_{1},y_{1},z_{1})
-l_{2}(t,x^{\prime}_{1},y^{\prime}_{1},x^{\prime\prime}_{1},y^{\prime\prime}_{2},x_{1},y_{2},z_{2})|,\\
&\leq L_{2}(|h_{2}(t,B_{2}(t)^{\top}y_{1}+\tau \bar{B}_{2}(t)^{\top}y^{\prime\prime}_{1}+D_{2}(t)^{\top}z_{1})\\
&\quad \quad\quad\quad \quad\quad\quad \quad\quad-h_{2}(t,B_{2}(t)^{\top}y_{2}+\tau \bar{B}_{2}(t)^{\top}y^{\prime\prime}_{2}+D_{2}(t)^{\top}z_{2}),\\
&l_{1}:=b_{1},\sigma_{1},f_{1}, l_{2}:=b_{2},\sigma_{2},f_{2}.
\end{cases}
\end{align}

\item[3)] Monotonicity conditions:  set
$V_{1}=(x_{1}, y_{1},z_{1}), V_{2}=(x_{2}, y_{2},z_{2}),\bar{V}_{1}=(\bar{x}_{1}, \bar{y}_{1},\bar{z}_{1}),\\ \bar{V}_{2}=(\bar{x}_{2}, \bar{y}_{2},\bar{z}_{2}),$ $\theta_{i}:=(x^{\prime}_{1},y^{\prime}_{1},x^{\prime}_{2},y^{\prime}_{2},x_{i},y_{i},z_{i}), $  $\bar{\theta}_{i}:=(\bar{x}^{\prime}_{1},\bar{y}^{\prime}_{1},\bar{x}^{\prime}_{2},\bar{y}^{\prime}_{2},\bar{x}_{i},\bar{y}_{i},\bar{z}_{i}),i=1,2. $  We assume that

\begin{align}
\begin{cases}
&\langle \Psi_{1}(y_{1},y)-\Psi_{1}(\bar{y}_{1},y), y_{1}-\bar{y}_{1} \rangle \leq -L_{3}|\bar{h}_{11}(\frac{H^{\top}y_{1}-\tau H^{\top}y}{1-\tau^{2}})-\bar{h}_{11}(\frac{H^{\top}\bar{y}_{1}-\tau H^{\top}y}{1-\tau^{2}} )|^{2}\\
&\quad \quad\quad\quad \quad\quad\quad \quad\quad\quad\quad\quad\quad\quad-L_{3}|\bar{h}_{12}(\frac{H^{\top}y-\tau H^{\top}y_{1}}{1-\tau^{2}})-\bar{h}_{12}(\frac{H^{\top}y-\tau H^{\top}\bar{y}_{1}}{1-\tau^{2}} )|^{2},\\
&\langle \Psi_{2}(y,y_{2})-\Psi_{2}(y,\bar{y}_{2}), y_{2}-\bar{y}_{2} \rangle
\leq -L_{3}|\bar{h}_{21}(\frac{H^{\top}y_{2}-\tau H^{\top}y}{1-\tau^{2}})-\bar{h}_{21}(\frac{H^{\top}\bar{y}_{2}-\tau H^{\top}y}{1-\tau^{2}} )|^{2}\\
&\quad \quad\quad\quad \quad\quad\quad \quad\quad\quad\quad\quad\quad\quad-L_{3}|\bar{h}_{22}(\frac{H^{\top}y-\tau H^{\top}y_{2}}{1-\tau^{2}})-\bar{h}_{22}(\frac{H^{\top}y-\tau H^{\top}\bar{y}_{2}}{1-\tau^{2}} )|^{2}\\
& \langle\Phi_{1}(x_{1},x_{2})-\Phi_{1}(\bar{x}_{1},x_{2}), x_{1}-\bar{x}_{1} \rangle
+\langle \Phi_{2}(x_{1},x_{2})-\Phi_{2}(x_{1},\bar{x}_{2}), x_{2}-\bar{x}_{2} \rangle \geq 0,\\
&\langle \Gamma_{1}(\theta_{1})-\Gamma_{1}(\bar{\theta}_{1}), V_{1}-\bar{V}_{1} \rangle+\langle \Gamma_{2}(\theta_{2})-\Gamma_{2}(\bar{\theta}_{2}), V_{2}-\bar{V}_{2} \rangle\\
& \leq -L_{3}|h_{1}(t,B_{1}(t)^{\top}y_{1}+\tau \bar{B}_{1}(t)^{\top}y^{\prime}_{1}+D_{1}(t)^{\top}z_{1})\\
&\quad\quad\quad\quad\quad\quad\quad\quad\quad-h_{1}(t,B_{1}(t)^{\top}\bar{y}_{1}+\tau \bar{B}_{1}(t)^{\top}\bar{y}^{\prime}_{1}+D_{1}(t)^{\top}\bar{z}_{1})|\\
&\quad-L_{3}|h_{2}(t,B_{2}(t)^{\top}y_{2}+\tau \bar{B}_{2}(t)^{\top}y^{\prime}_{2}+D_{2}(t)^{\top}z_{2})\\
&\quad\quad\quad\quad\quad\quad\quad\quad\quad-h_{2}(t,B_{2}(t)^{\top}\bar{y}_{2}+\tau \bar{B}_{2}(t)^{\top}\bar{y}^{\prime}_{2}+D_{2}(t)^{\top}\bar{z}_{2})|.
\end{cases}
\end{align}

\end{itemize}

\end{itemize}
Finally, we make the following convention:
 The letter $C$  will denote  positive constant which only depends on the constants in Assumption, whose value may
vary from one place to another.
\subsection{Existence and uniqueness  of extend MF-FBSDEs}

The main results of this section  will be summarized as follows and the proof of  the following  results will be listed   in Appendix.
\subsection{Prior estimates}
\bl\label{t1}
Assume  that the the coefficients satisfy $\mathbf{Assumption 1}.$  We also assume that \\ $V_{i}(\cdot) \in \mM^{2}_{\mF}(\mR^{{n+n+nd}}), i=1,2 $  is a solution to System $(\pi),$ where $V_{i}(\cdot):=(X_{i}(t)^{\top}(\cdot),$ $Y_{i}(\cdot)^{\top},$\\$Z_{i}(\cdot)^{\top})^{\top}, i=1,2.$  Consequently, we obtain the following estimate:
\begin{equation}\label{11}
\mE[\Lambda_{V_{1}}]+\mE[\Lambda_{V_{2}}]\leq C(\Psi(0,0)+\Psi(0,0)+|\Phi_{1}(0,0)|^{2}+|\Phi_{2}(0,0)|^{2}+\Xi_{\Gamma_{1}(\cdot,0)}+\Xi_{\Gamma_{2}(\cdot,0)}),
\end{equation}
where $\Xi_{\Gamma_{i}(\cdot,0)}, \Lambda_{V_{i}}, i=1,2$  are given by equations \eqref{99} and \eqref{301}, respectively, $C$ is a constant only depending on  the constants in $\mathbf{Assumption 1}.$ Consider another set of coefficients $(\tilde{\Psi}_{i},\tilde{\Phi}_{i}, \tilde{\Gamma}_{i}),i=1,2$ and let $\tilde{V}_{i}(\cdot):=(\tilde{X}_{i}(t)^{\top}(\cdot),\tilde{Y}_{i}(\cdot)^{\top},\tilde{Z}_{i}(\cdot)^{\top})^{\top}   \in \mM^{n+n+nd}, i=1,2 $ be a solution to System $(\pi)$ associated with $(\tilde{\Psi}_{i}, \tilde{\Phi}_{i}, \tilde{\Gamma}_{i}).$ Then the following estimate holds:
\begin{align}\label{12}
\mE[\Lambda_{\hat{\theta}_{1}}]+\mE[\Lambda_{\hat{\theta}_{2}}]\leq C(|\hat{\Psi}(\tilde{Y}_{1}(0),\tilde{Y}_{2}(0))|^{2}+\mE(|\hat{\Phi}(\tilde{X}_{1}(T),\tilde{X}_{2}(T))|^{2}+\Xi_{\hat{\Gamma}_{1}(\cdot, \tilde{\theta}_{1}(\cdot))})+\Xi_{\hat{\Gamma}_{2}(\cdot, \tilde{\theta}_{2}(\cdot))})),
\end{align}
where $\hat{l}=l-\tilde{l}, l=\theta_{i},\Psi_{i},\Phi_{i},\Gamma_{i},$ $$\tilde{\theta}_{i}(t):=((\mE[\tilde{X}_{1}(t)])^{\top}, (\mE[\tilde{Y}_{1}(t)])^{\top}, (\mE[\tilde{X}_{2}(t)])^{\top}, (\mE[\tilde{Y}_{2}(t)])^{\top},  \tilde{ X}_{i}(t)^{\top}, \tilde{Y}_{i}(t)^{\top}, \tilde{Z}_{i}(t)^{\top})^{\top}, i=1,2. $$
\el

\subsection{ Method of continuity}  In this part, we will introduce the method of continuity  and give a useful lemma which will be used to prove the existence and uniqueness result.
Under the Assumption 1, we define another set of coefficients $(\Psi^{0}_{i}, \Phi^{0}_{i}, \Gamma^{0}_{i}), i=1,2$  with

\begin{align}\label{21}
&\Psi^{0}_{1}(y_{1},y_{2}):=H\bar{h}_{11}(\frac{H^{\top}y_{1}-\tau H^{\top}y_{2}}{1-\tau^{2}})+H\bar{h}_{12}(\frac{H^{\top}y_{2}-\tau H^{\top}y_{1}}{1-\tau^{2}}),\no\\
&\Psi^{0}_{2}(y_{1},y_{2}):=H\bar{h}_{21}(\frac{H^{\top}y_{2}-\tau H^{\top}y_{1}}{1-\tau^{2}})+H\bar{h}_{22}(\frac{H^{\top}y_{1}-\tau H^{\top}y_{2}}{1-\tau^{2}}),\no\\
&  \Phi^{0}_{i}(x_{1},x_{2}):=0, f^{0}_{i}(t,\theta_{i}):=0,\no\\
&b^{0}_{i}(t,\theta_{i}):=B_{i}(t)h_{i}(t, B_{i}(t)^{\top}y_{i}+\tau\bar{B}_{i}(t)^{\top}y'_{i}+D_{i}(t)^{\top}z_{i}),\no\\
 &\sigma^{0}_{i}(t,\theta_{i}):=D_{i}(t)h_{i}(t, B_{i}(t)^{\top}y_{i}+\tau\bar{B}_{i}(t)^{\top}y'_{i}+D_{i}(t)^{\top}z_{i}),i=1,2,
\end{align}
for any $\theta_{i}:=(x_{1}^{\prime\top},y_{1}^{\prime\top}, x_{2}^{\prime\top}, y_{2}^{\prime\top}, x_{i}^{\top},y_{i}^{\top},z_{i}^{\top})^{\top}\in \mR^{n+n+n+n+n+n+nd}, i=1,2.$  We employ the notations $\Gamma^{0}_{i}:=((f^{0}_{i})^{\top},(b^{0}_{i})^{\top}, (\sigma^{0}_{i})^{\top})^{\top} $ and $ \sigma^{0}_{i}:=((\sigma^{0}_{i1})^{\top},(\sigma^{0}_{i2})^{\top},\cdots, (\sigma^{0}_{id})^{\top})^{\top}$ in this context. A simple verification demonstrates that the coefficients $(\Psi^{0}_{i}, \Phi^{0}_{i}, \Lambda^{0}_{i})$  adhere to Assumption 1 with  the same parameters(If necessary, make appropriate adjustments, but ensure they do not affect the entire proof process).

For any $\xi_{i}\in \mR^{n},  \zeta_{i} \in L^{2}_{\sF_{T}}(\mR^{d}),$ and any $\beta_{i}(\cdot):=(\phi_{i}(\cdot)^{\top}, \psi_{i}(\cdot)^{\top},\gamma_{i}(\cdot)^{\top})^{\top}\in \sM_{\mF}^{2}(\mR^{n+n+nd})$ and $\gamma_{i}:=((\gamma_{i1})^{\top},(\gamma_{i2})^{\top},\cdots, (\gamma_{id})^{\top})^{\top},i=1,2,$ we now proceed to introduce a family of MF-FBSDEs defined by the parameter $\alpha\in [0,1],$
\begin{align}\label{23}
\begin{cases}
&\dif X_{1}^{\alpha}(t)=[b_{1}^{\alpha}(t,\theta_{1}^{\alpha}(t))+\psi_{1}(t)]\dif t
+\sum^{d}_{j=1}[\sigma_{1j}^{\alpha}(t,\theta^{\alpha}_{1}(t))+\gamma_{1j}(t)]\dif W_{j}(t), t\in [0,T]\\
&\dif Y^{\alpha}_{1}(t)^{\top}(t) =\{[f^{\alpha}_{1}(t,\theta_{1}^{\alpha}(t))+\phi_{1}(t)]\dif t+\sum^{d}_{j=1} Z^{\alpha}_{1j}(t)\dif W_{j}(t)\\
& X_{1}^{\alpha}(0)=\Psi_{1}^{\alpha}(Y_{1}^{\alpha}(0),Y_{2}^{\alpha}(0))+ \xi_{1},  Y^{\alpha}_{1}(T)=\Phi_{1}^{\alpha}(X^{\alpha}_{1}(T),X^{\alpha}_{2}(T) )+\zeta_{1},
\end{cases}
\end{align}

 \begin{align}\label{24}
\begin{cases}
&\dif X_{2}^{\alpha}(t)=[b_{2}^{\alpha}(t,\theta_{2}^{\alpha}(t))+\psi_{2}(t)]\dif t
 +\sum^{d}_{j=1}[\sigma_{2j}(t,\theta^{\alpha}_{2}(t))+\gamma_{2j}(t)]\dif W_{j}(t), t\in [0,T]\\
&\dif Y_{2 }^{\alpha}(t) =[f_{2}^{\alpha}(t,\theta_{2}^{\alpha}(t))+\phi_{2}(t)]\dif t+\sum^{d}_{j=1} Z^{\alpha}_{2j}(t)\dif W_{j}(t)\\
& X_{2 }^{\alpha}(0)=\Psi_{2}^{\alpha}(Y_{1}^{\alpha}(0),Y_{2}^{\alpha}(0))+ \xi_{2},  Y_{2}^{\alpha}(T)=\Phi_{2}^{\alpha}(X_{1}^{\alpha}(T),X_{2}^{\alpha}(T) )+\zeta_{2},
\end{cases}
\end{align}
 where $(\Psi_{i}^{\alpha}, \Phi_{i}^{\alpha}, \Gamma_{i}^{\alpha}):=\alpha(\Psi_{i}, \Phi_{i}, \Gamma_{i})+(1-\alpha)(\Psi_{i}^{0}, \Phi_{i}^{0}, \Gamma_{i}^{0}), i=1,2,$\\   $\theta^{\alpha}_{i}(t):=((\mE[X^{\alpha}_{1}(t)])^{\top}, (\mE[Y^{\alpha}_{1}(t)])^{\top}, (\mE[X^{\alpha}_{2}(t)])^{\top}, (\mE[Y^{\alpha}_{2}(t)])^{\top},   X^{\alpha}_{i}(t)^{\top}, Y^{\alpha}_{i}(t)^{\top}, Z^{\alpha}_{i}(t)^{\top})^{\top}. $ Obviously, coefficients $(\Psi_{i}^{\alpha}+\xi_{i}, \Phi_{i}^{\alpha}+\zeta_{i}, \Gamma_{i}^{\alpha}+ \beta_{i}), i=1,2$ of Eq.\eqref{24} satisfy $\mathrm{Assumption\, 1}.$
For simplicity,  Eq.\eqref{23} and Eq.\eqref{24} are denoted by System $(\pi_{1}).$

We see two extreme cases. For the case $\alpha = 1,$  $(\xi_{i}, \zeta_{i}, \beta_{i}(\cdot)), i=1,2$ are all vanish,  System $(\pi_{1})$  reduces to System $(\pi)$, which  we aim to investigate.  For the case $\alpha = 0,$  system $(\pi_{1})$ degenerates into a decoupled form, which can be solved using the established results for SDEs (see \cite{Zhang2017}) and BSDEs (see \cite{PardouxPeng1990}).

\bl\label{t4}
Assume that Assumption 1 is satisfied for the coefficients  $(\Psi_{i}, \Phi_{i}, \Gamma_{i}),i=1,2.$ We can find an absolute constant $\delta_{0} > 0$ such that if  for some $\alpha_{0}\in [0,1),$ System $(\pi_{1})$ is uniquely solvable in $\mM_{\mF}^{2}(\mR^{n+n+nd})\times \mM_{\mF}^{2}(\mR^{n+n+nd}) $ for any $(\xi_{i}, \zeta_{i}, \beta_{i}(\cdot)) \in \mR^{n}\times L^{2}_{\sF_{T}}(\mR^{n})\times \sM_{\mF}^{2}(\mR^{n+n+nd}), i=1,2,$  then replacing $\alpha_0$ by any $\alpha \in (\alpha_0, (\alpha_0 + \delta_0) \wedge
1]$, the same conclusion remains true.

\el

The well-posed of System $(\pi)$ be listed as follows.
\bt\label{tt}
Provided that Assumption 1 holds for the coefficients $(\Psi_{i},\Phi_{i},\Gamma_{i}), i=1,2 ,$ System $(\pi)$ admits a unique solution.
\et

\section{Application to stochastic linear-convex problems}
Here, we investigate $Problem (LC),$ a problem formulated in Section 1 based on the solvability of extended  MF-FBSDEs. First, we state the rigorous assumptions required for the coefficients in the controlled system \eqref{4}.\\
$\mathbf{Assumption\, 2:}$\\
 $A_{i}(\cdot), \bar{A}_{i}(\cdot), C_{ij}(\cdot) \in L^{\infty}(\mathbb{R}^{n \times n}), B_{i}(\cdot), \bar{B}_{i}(\cdot), D_{ij}(\cdot) \in L^{\infty}(\mathbb{R}^{n \times k}), H \in \mathbb{R}^{n \times m}, \rho_{i}(\cdot) \in L_{\mathbb{F}}^{2}(\mathbb{R}^n),$ $\kappa_{ij}(\cdot) \in L_{\mathbb{F}}^{2}(\mathbb{R}^n)$, and $x_0 \in \mathbb{R}^n$ for any $j = 1, 2, \dots, d, i=1,2$.  Moreover, there exists a sufficiently small constant $\tau_{1}> 0$ such that $\sup_{t\in [0,T]}|\bar{A}_{i}(t)|<\tau_{1},i=1,2. $

For convenience, we introduce the notations $C_{i}(\cdot) := (C_{i1}(\cdot)^{\top}, C_{i2}(\cdot)^{\top}, \dots, C_{id}(\cdot)^{\top})^{\top}, D_{i}(\cdot)$ $:= (D_{i1}(\cdot)^{\top}, D_{i2}(\cdot)^{\top}, \dots, D_{id}(\cdot)^{\top})^{\top}$, and $\kappa_{i}(\cdot) := (\kappa_{i1}(\cdot)^{\top}, \kappa_{i1}(\cdot)^{\top}, \dots, \kappa_{id}(\cdot)^{\top})^{\top}$.  For any $(\xi_{i}, u_{i}(\cdot)) \in \mathbb{R}^m \times L_{\mathbb{F}}^{2}(\mathbb{R}^k)$,   as a special case of Theorem \ref{tt}, there exists a unique solution to Eq.\eqref{4} $X_{i}(\cdot) \equiv X_{i}(\cdot; \xi_{i}, u_{i}(\cdot)) \in S_{\mathbb{F}}^{2}(\mathbb{R}^n),i=1,2.$

\subsection{Convex Criterion Functional}
For a function $g:\mR^{n}\supset D\rightarrow \mR,$    we denote that
$$\nabla g(x) := \left( \frac{\partial g}{\partial x_1}(x), \frac{\partial g}{\partial x_2}(x), \dots, \frac{\partial g}{\partial x_n}(x) \right)^\top.$$

$\mathbf{Assumption\, 3:}$
\begin{itemize}
\item[1)]  $f_{2i}(\cdot)$ and $f_{3i}( t,\cdot),i=1,2$ are convex, and $f_{1i}(\cdot)$ and $f_{4i}(t,\cdot),i=1,2$ are uniformly convex with parameter $\delta > 0$ (see Definition \ref{d1} in the Appendix) for almost all $(\omega, t)\in \Omega \times [0,T]$.
  \item[2)] $f_{1i}(\cdot)$, $f_{2i}(\cdot)$, $f_{3i}(t, \cdot)$ and $f_{4i}(t, \cdot)$ are continuously differentiable for almost all $(\omega, t)\in \Omega \times [0,T]$. Moreover, $\nabla f_{2i}(\cdot),$  $\nabla f_{3i}(t,\cdot), i=1,2$ are uniform Lipschitz continuous in $(\omega, t)$ and the Lipschitz constants of   $\nabla f_{2i}(\cdot)$ are small enough.
 \item[3)] $f_{2i}$, $f_{3i}$, and $f_{4i}$ are $\sF_{T}\times\sB(\mathbb{R}^n)$-measurable, $\sP \times \mathcal{B}(\mathbb{R}^n)$-measurable, and $\sP \times \sB(\mathbb{R}^k)$-measurable, respectively.
\item[4)]$f_{2i}(0)\in L_{\sF_{T}}(\mR)$ $f_{3i}(\cdot,0),f_{4i}(\cdot,0)\in L_{\mF}(\mathbb{R}), \nabla f_{2i}(0) \in L_{\mathcal{F}_T}^2(\mathbb{R}^n), \nabla f_{3i}(\cdot,0) \in L_{\mathbb{F}}^{2}(\mathbb{R}^n), \text{and} \\ \nabla f_{4i}(\cdot,0) \in L_{\mathbb{F}}^2(\mathbb{R}^k),i=1,2.$

\end{itemize}

As we can see under Assumption 3, Lemma \ref{A1} implies that both $\nabla f_{1i}(\cdot): \mathbb{R}^{m} \rightarrow \mathbb{R}^{m}, i=1,2$ and $\nabla f_{4i}(\cdot): \mathbb{R}^{k} \rightarrow \mathbb{R}^{k},i=1,2$ are bijective. Let $(\nabla f_{1i})^{-1}(\cdot), i=1,2$ and $(\nabla f_{4i})^{-1}( \cdot), i=1,2$ denote the four inverse mappings, respectively. Moreover, under the Assumption 3, we also know that
\begin{align}\label{61}
\mathbb{E}&\bigg\{\bigg.f_{21}(X_{1}(T)+X_{2}(T) ) +f_{22}( X_{1}(T)+X_{2}(T) )\no\\
  &\quad\quad\quad\quad\quad\quad\quad\quad\quad+ \int^{T}_{0}f_{31}(t,X_{1}(t) )\dif t +\int^{T}_{0}f_{32}(t,X_{2}(t))\dif t\bigg\}\bigg. < \infty.
\end{align}
for any $X_{1}(\cdot), X_{2}(\cdot) \in S_{\mF}^{2}(\mathbb{R}^n)$. However, the final integral in the criterion functional \eqref{5} maybe diverge to $\infty$. To facilitate later analysis, we  present
\begin{align}\label{62}
\mathcal{U}_{1} := \left\{u_{1}(\cdot) \in L_{\mF}^2(\mathbb{R}^k) \ \bigg| \ \mathbb{E} \int_0^T f_{41}(t, u_{1}(t)) dt < \infty\right\}. \
\end{align}

\begin{align}\label{62+}
\mathcal{U}_{2} := \left\{u_{2}(\cdot) \in L_{\mF}^2(\mathbb{R}^k) \ \bigg| \ \mathbb{E} \int_0^T f_{42}(t, u_{2}(t)) dt < \infty \right\}.
\end{align}
It's apparent that, given $u_{1}(\cdot), u_{2}(\cdot) \in L_{\mF}^2(\mathbb{R}^k)$, then $0 \leq J(\xi_{1},\xi_{2}, u_{1}(\cdot),u_{2}(\cdot) ) < \infty$ if and only if $(u_{1}(\cdot), u_{2}(\cdot)) \in \mathcal{U}_{1}\times \mathcal{U}_{2} $.  It's evident that the set  $\mathcal{U}_{1}, \mathcal{U}_{2}$ defined by \eqref{62}, \eqref{62+} are nonempty. We define $\mathbb{R}^m\times \mathbb{R}^m \times \mathcal{U}_{1}\times \mathcal{U}_{2}$ as the admissible control set. When $(\xi_{1}, \xi_{2},u_{1}(\cdot), u_{2}(\cdot)) \in \mathbb{R}^m \times \mathbb{R}^m \times \mathcal{U}_{1}\times \mathcal{U}_{2}$,  we refer to  them  as quadruple admissible controls. Moreover, $X_{i}(\cdot) \equiv X_{i}(\cdot; \xi_{i}, u_{i}(\cdot)),i=1,2$ and $(\xi_{1}, u_{1}(\cdot), X_{1}(\cdot),\xi_{2}, u_{2}(\cdot), X_{2}(\cdot))$ are called the corresponding admissible state and an admissible sextet, respectively.

Unlike most stochastic optimal control problems, Problem (LC) in this article includes two  initial controls $\xi_{1},\xi_{2}$ and two process controls $u_{1},u_{2}$.   Nevertheless, it remains a form of Bolza problem.

\subsection{Stochastic Hamiltonian system}The optimal control quartet of Problem (LC) will be characterized using an extended MF-FBSDE, also referred to as a stochastic Hamiltonian system in control theory.

\begin{equation}\label{22++}
J(\xi_{1},\xi_{2}, u_{1}(\cdot), u_{2}(\cdot)) - J(\xi^*_{1},\xi^*_{2}, u^*_{1}(\cdot),u^*_{2}(\cdot)) = \sum_{i=1}^{4} \Diamond_i
\end{equation}

where
\begin{align*}
\Diamond_1 &:= f_{11}(\xi_{1}+\tau\xi_{2} ) - f_{11}(\xi^{*}_{1}+\tau\xi^{*}_{2})+f_{12}(\xi_{2}+\tau\xi_{1} ) - f_{12}(\xi^{*}_{2}+\tau\xi^{*}_{1}), \\
\Diamond_2 &:= \sum^{2}_{i=1}\{\mathbb{E} \left[ f_{2i}(X_{1}(T)+ X_{2}(T)) - f_{2i}(X^*_{1}(T)+X^*_{2}(T) ) \right]\}, \\
\Diamond_3 &:= \sum^{2}_{i=1}\bigg\{\bigg.\mathbb{E} \int_{0}^{T} \left[ f_{3i}(X_{i}(t) ) - f_{3i}(X^*_{i}(t) ) \right] dt\bigg\}\bigg., \\
\Diamond_4 &:=  \sum^{2}_{i=1}\bigg\{\bigg.\mathbb{E} \int_{0}^{T} \left[ f_{4i}(t, u_{i}(t)) - f_{4i}(t, u^*_{i}(t)) \right] dt\bigg\}\bigg..
\end{align*}

Due to the convexity of $f_{2i}$ and $f_{3i}$ and the uniform convexity of $f_{1i}$ and $f_{4i}$, Lemma \ref{100} (2) in the Appendix works to yield

\begin{align}\label{27}
&\Diamond_1 \geq \langle \nabla f_{11}(\xi^*_{1}+\tau\xi^*_{2}), \xi_{1} - \xi^*_{1} \rangle + \langle \tau\nabla f_{11}(\xi^*_{1}+\tau\xi^*_{2}), \xi_{2} - \xi^*_{2} \rangle\no\\
&\quad\quad+\langle \nabla f_{12}(\tau\xi^*_{1}+\xi^*_{2}), \xi_{2} - \xi^*_{2} \rangle + \langle \tau\nabla f_{12}(\tau\xi^*_{1}+\xi^*_{2}), \xi_{1} - \xi^*_{1} \rangle\no\\
&\quad\quad+\frac{\delta}{2} |\xi_{1}+\tau\xi_{2} - \xi^*_{1}-\tau\xi^*_{2}|^2+\frac{\delta}{2} |\tau\xi_{1}+\xi_{2} - \tau\xi^*_{1}-\xi^*_{2}|^2,\no\\
&\Diamond_2\geq \sum^{2}_{i=1}\bigg\{\bigg.\mathbb{E} [ \langle \nabla f_{2i}(X^*_{1}(T)+X^*_{2}(T) ), (X_{1}(T)+X_{2}(T)) - (X^*_{1} (T)+(X^*_{2} (T)\rangle]\bigg\}\bigg.,\no \\
 &\Diamond_3
 \geq \sum^{2}_{i=1}\mathbb{E} \int_{0}^{T} \left[ \langle \nabla f_{3i}( t,X^*_{i}(t)), X_{i}(t) - X^*_{i}(t) \rangle \right] dt \no\\
&\Diamond_4\geq \sum^{2}_{i=1}\mathbb{E} \int_{0}^{T} \left[ \langle \nabla f_{4i}(t, u^*_{i}(t)), u_{i}(t) - u^*_{i}(t) \rangle + \frac{\delta}{2} |u_{i}(t) - u^*_{i}(t)|^2 \right] dt.
\end{align}
Clearly, \eqref{27} restates \eqref{22++} as an inequality. To make this inequality more tractable, we employ a duality viewpoint and introduce a BSDE. This BSDE, defined on $[0, T],$ takes the following form:
\begin{align}\label{28}
&\dif Y_{1}(t)=-[\nabla f_{31}(t,X^{*}_{1}(t))+A_{1}^{\top}(t)Y_{1}(t)
+C_{1}^{\top}(t)Z_{1}(t)+\bar{A}_{1}^{\top}(t)\mE[Y_{1}(t)]+\bar{A}_{2}^{\top}(t)\mE[Y_{2}(t)] ]\dif t\no\\
&\quad\quad\quad+\sum^{d}_{j=1}Z_{1j}(t)\dif W_{j}(t), t\in [0,T],\no\\
& Y_{1}(T)=\nabla f_{21}(X^{*}_{1}(T)+X^{*}_{2}(T))+\nabla f_{22}(X^{*}_{1}(T)+X^{*}_{2}(T)),
\end{align}
and

\begin{align}\label{29}
&\dif Y_{2}(t)=-[\nabla f_{32}(t,X^{*}_{2}(t))+A_{2}^{\top}(t)Y_{2}(t)
+C_{2}^{\top}(t)Z_{2}(t)+\bar{A}_{2}^{\top}(t)\mE[Y_{1}(t)] ]\dif t\no\\
&\quad\quad\quad+\sum^{d}_{j=1}Z_{2j}(t)\dif W_{j}(t), t\in [0,T],\no\\
&  Y_{2}(T)=\nabla f_{21}(X^{*}_{1}(T)+X^{*}_{2}(T))+\nabla f_{22}(X^{*}_{1}(T)+X^{*}_{2}(T)).
\end{align}
Eq.\eqref{28} has a unique solution $(Y_{1},Z_{1})\in L^{2}_{\mF}(\mR^{n})\times L^{2}_{\mF}(\mR^{nd})$ and Eq.\eqref{29} has a unique solution $(Y_{2},Z_{2})\in L^{2}_{\mF}(\mR^{n})\times L^{2}_{\mF}(\mR^{nd}).$
Applying It\^{o}'s formula to $\langle Y_{1}(t), X_{1}(t)-X^{*}_{1}(t) \rangle, $   we have
\begin{align}\label{30}
&\mE\bigg\{\bigg. \langle \nabla f_{21}(X^{*}_{1}(T)+X^{*}_{2}(T)),   X_{1}(T)-X^{*}_{1}(T)\rangle+\langle \nabla f_{22}(X^{*}_{1}(T)+X^{*}_{2}(T)),   X_{1}(T)-X^{*}_{1}(T)\rangle  \no\\
  &\quad \quad \quad\quad \quad \quad\quad \quad \quad+\int^{T}_{0}\langle \nabla f_{31}(X^{*}_{1}(t)),   X_{1}(t)-X^{*}_{1}(t)\rangle \dif t \bigg\}\bigg.\no\\
&=\langle H^{\top}Y_{1}(0), \xi_{1}-\xi^{*}_{1}\rangle-\mE\int^{T}_{0}\langle A_{1}(t)^{\top}Y_{1}(t), X_{1}(t)-X^{*}_{1}(t) \rangle\dif t\no\\
&-\mE\int^{T}_{0}\langle \bar{A}_{1}(t)^{\top}\mE[Y_{1}(t)], X_{1}(t)-X^{*}_{1}(t) \rangle\dif t-\mE\int^{T}_{0}\langle \bar{A}_{2}(t)^{\top}\mE[Y_{2}(t)], X_{1}(t)-X^{*}_{1}(t) \rangle\dif t\no\\
&+\mE\int^{T}_{0}\langle A_{1}(t)^{\top}Y_{1}(t), X_{1}(t)-X^{*}_{1}(t) \rangle\dif t+\mE\int^{T}_{0}\langle \bar{A}_{2}(t)^{\top}\mE[Y_{1}(t)], X_{2}(t)-X^{*}_{2}(t) \rangle\dif t\no\\
&+\mE\int^{T}_{0}\langle \bar{A}^{\top}_{1}(t)\mE[Y_{1}(t)], X_{1}(t)-X^{*}_{1}(t) \rangle\dif t\no\\
&+\mE\int^{T}_{0}\langle B_{1}^{\top}(t)Y_{1}(t)+\tau \bar{B}_{1}^{\top}(t)\mE[Y_{1}(t)]+ D_{1}^{\top}(t)Z_{1}(t), u_{1}(t)-u^{*}_{1}(t)   \rangle\dif t.
\end{align}
Similarly, applying It\^{o}'s formula to $\langle Y_{2}(t), X_{2}(t)-X^{*}_{2}(t) \rangle$  it yields

\begin{align}\label{31}
&\mE\bigg\{\bigg. \langle \nabla f_{21}(X^{*}_{1}(T)+X^{*}_{2}(T)),   X_{2}(T)-X^{*}_{2}(T)\rangle+\langle \nabla f_{22}(X^{*}_{1}(T)+X^{*}_{2}(T)),   X_{2}(T)-X^{*}_{2}(T)\rangle  \no\\
  &\quad \quad \quad\quad \quad \quad\quad \quad \quad+\int^{T}_{0}\langle \nabla f_{32}(X^{*}_{2}(t)),   X_{2}(t)-X^{*}_{2}(t)\rangle \dif t \bigg\}\bigg.\no\\
&=\langle H^{\top}Y_{2}(0), \xi_{2}-\xi^{*}_{2}\rangle-\mE\int^{T}_{0}\langle A_{2}(t)^{\top}Y_{2}(t), X_{2}(t)-X^{*}_{2}(t) \rangle\dif t\no\\
&-\mE\int^{T}_{0}\langle \bar{A}_{2}(t)^{\top}\mE[Y_{1}(t)], X_{2}(t)-X^{*}_{2}(t) \rangle\dif t \rangle\dif t\no\\
&+\mE\int^{T}_{0}\langle A_{2}(t)^{\top}Y_{2}(t), X_{2}(t)-X^{*}_{2}(t) \rangle\dif t+\mE\int^{T}_{0}\langle \bar{A}_{2}(t)^{\top}\mE[Y_{2}(t)], X_{1}(t)-X^{*}_{1}(t) \rangle\dif t\no\\
& +\mE\int^{T}_{0}\langle    B_{2}^{\top}(t)Y_{2}(t)+ \tau\bar{B}_{2}^{\top}(t)\mE[Y_{2}(t)]+D_{2}^{\top}(t)Z_{2}(t),   u_{2}(t)-u^{*}_{2}(t) \rangle\dif t.
\end{align}
Summing up equations \eqref{30}-\eqref{31} yields
\begin{align}\label{34}
&\mE\bigg\{\bigg. \nabla f_{21}(X^{*}_{1}(T)+X^{*}_{2}(T)),   (X_{1}(T)+X_{2}(T))-(X^{*}_{1}(T)+X^{*}_{2}(T))\rangle\no\\
& \quad\quad +\nabla f_{22}(X^{*}_{1}(T)+X^{*}_{2}(T)),     (X_{1}(T)+X_{2}(T))-(X^{*}_{1}(T)+X^{*}_{2}(T))\rangle\no\\
&\quad\quad+\int^{T}_{0}\langle \nabla f_{31}(X^{*}_{1}(t)),   X_{1}(t)-X^{*}_{1}(t)\rangle \dif t +\int^{T}_{0}\langle \nabla f_{32}(X^{*}_{2}(t)),   X_{2}(t)-X^{*}_{2}(t)\rangle \dif t\bigg\}\bigg.\no\\
  &=\langle H^{\top}Y_{1}(0), \xi_{1}-\xi^{*}_{1}\rangle
+\langle H^{\top}Y_{2}(0), \xi_{2}-\xi^{*}_{2}\rangle\no\\
&+\mE\int^{T}_{0}\langle B_{1}^{\top}(t)Y_{1}(t)+\tau \bar{B}_{1}^{\top}(t)\mE[Y_{1}(t)]+ D_{1}^{\top}(t)Z_{1}(t), u_{1}(t)-u^{*}_{1}(t)   \rangle\dif t\no\\
& +\mE\int^{T}_{0}\langle    B_{2}^{\top}(t)Y_{2}(t)+ \tau\bar{B}_{2}^{\top}(t)\mE[Y_{2}(t)]+D_{2}^{\top}(t)Z_{2}(t),   u_{2}(t)-u^{*}_{2}(t) \rangle\dif t.
\end{align}
Combining \eqref{22++}-\eqref{34} leads to
\begin{align}\label{35}
&J(\xi_{1},\xi_{2}, u_{1}(\cdot), u_{2}(\cdot)) - J(\xi^*_{1},\xi^*_{2}, u^*_{1}(\cdot),u^*_{2}(\cdot))\no\\
&\geq\langle \nabla f_{11}(\xi_{1}^{*}+\tau\xi_{2}^{*})+\tau \nabla f_{12}(\tau\xi_{1}^{*}+\xi_{2}^{*})+H^{\top}Y_{1}(0), \xi_{1}-\xi^{*}_{1}\rangle\no\\
&\quad +\langle\tau\nabla f_{11}(\xi_{1}^{*}+\tau\xi_{2}^{*})+\nabla f_{12}(\tau\xi_{1}^{*}+\xi_{2}^{*})+H^{\top}Y_{2}(0), \xi_{2}-\xi^{*}_{2}\rangle\no\\
&\quad +\mE\int^{T}_{0}\langle \nabla f_{41}(u^{*}_{1}(t))+ B^{\top}_{1}(t)Y_{1}(t)+ \tau\bar{B}_{1}^{\top}(t)\mE[Y_{1}(t)]+D^{\top}_{1}(t)Z_{1}(t), u_{1}(t)-u^{*}_{1}(t)   \rangle\dif t\no\\
&\quad +\mE\int^{T}_{0}\langle \nabla f_{42}(u^{*}_{2}(t))+B^{\top}_{2}(t)Y_{2}(t)+ \tau\bar{B}_{2}^{\top}(t)\mE[Y_{2}(t)]+D^{\top}_{2}(t)Z_{2}(t),   u_{2}(t)-u^{*}_{2}(t) \rangle\dif t\no\\
&\quad+\frac{\delta}{2}\bigg\{\bigg. |\xi_{1}+\tau\xi_{2} - \xi^*_{1}-\tau\xi^*_{2}|^2+ |\tau\xi_{1}+\xi_{2} - \tau\xi^*_{1}-\xi^*_{2}|^2\no\\
& \quad+\mE\int^{T}_{0} |u_{1}(t)-u^{*}_{1}(t)|^{2}\dif t +\mE\int^{T}_{0} |u_{2}(t)-u^{*}_{2}(t)|^{2}\dif t    \bigg\}\bigg.
\end{align}
Denote $$\Diamond^{*}_7:=\nabla f_{11}(\xi_{1}^{*}+\tau\xi_{2}^{*})+\tau \nabla f_{12}(\tau\xi_{1}^{*}+\xi_{2}^{*})+H^{\top}Y_{1}(0),$$
$$\Diamond^{*}_8:=\tau\nabla f_{11}(\xi_{1}^{*}+\tau\xi_{2}^{*})+\nabla f_{12}(\tau\xi_{1}^{*}+\xi_{2}^{*})+H^{\top}Y_{2}(0),$$
$$\Diamond^{*}_9:= \nabla f_{41}(u^{*}_{1}(t))+ B^{\top}_{1}(t)Y_{1}(t)+ \tau\bar{B}_{1}^{\top}(t)\mE[Y_{1}(t)]+D^{\top}_{1}(t)Z_{1}(t),$$
$$\Diamond^{*}_{10}:=\nabla f_{42}(u^{*}_{2}(t))+B^{\top}_{2}(t)Y_{2}(t)+ \tau\bar{B}_{2}^{\top}(t)\mE[Y_{2}(t)]+D^{\top}_{2}(t)Z_{2}(t).$$
This section focuses on identifying optimal control pairs for Problem (LC). Considering inequality \eqref{35}, the optimality of $(\xi_{1},\xi_{2},u_{1},u_{2})$ appears equivalent to the following condition
$$(\Diamond^{*}_7, \Diamond^{*}_8, \Diamond^{*}_9, \Diamond^{*}_{10}) = (0, 0,0,0).$$
Then,  we have
\begin{align}\label{36}
\begin{cases}
\xi^{*}_{1}=\frac{1}{1-\tau^{2}}\bigg\{\bigg.(\nabla f_{11})^{-1}(-\frac{H^{\top}Y_{1}(0)-\tau H^{\top}Y_{2}(0)}{1-\tau^{2}})-\tau\nabla (f_{12})^{-1}(-\frac{H^{\top}Y_{2}(0)-\tau H^{\top}Y_{1}(0)}{1-\tau^{2}})\bigg\}\bigg.,\\
\xi^{*}_{2}=\frac{1}{1-\tau^{2}}\bigg\{\bigg.(\nabla f_{12})^{-1}(-\frac{H^{\top}Y_{2}(0)-\tau H^{\top}Y_{1}(0)}{1-\tau^{2}})-\tau\nabla (f_{11})^{-1}(-\frac{H^{\top}Y_{1}(0)-\tau H^{\top}Y_{2}(0)}{1-\tau^{2}})\bigg\}\bigg.,\\
u^{*}_{1}(\cdot)=(\nabla f_{41})^{-1}(-( B^{\top}_{1}(\cdot)Y_{1}(\cdot)+ \tau\bar{B}_{1}^{\top}(\cdot)\mE[Y_{1}(\cdot)]+D^{\top}_{1}(t)Z_{1}(\cdot))),\\
u^{*}_{2}(\cdot)=(\nabla f_{42})^{-1}(-( B^{\top}_{2}(\cdot)Y_{2}(\cdot)+ \tau\bar{B}_{2}^{\top}(\cdot)\mE[Y_{2}(\cdot)]+D^{\top}_{2}(\cdot)Z_{2}(\cdot))).
\end{cases}
\end{align}
This conjecture will be formally demonstrated in the following main result in this section.

At the conclusion of this section, we combine the SDEs for the state $X^{*}_{1}(\cdot), X^{*}_{2}(\cdot)$ [see \eqref{4}], the BSDEs \eqref{28}, \eqref{29} and the expression \eqref{36} into a single system, referred to as a stochastic Hamiltonian system, as follows (where the argument $t$ is omitted for simplicity):
\begin{align}\label{37}
\begin{cases}
&\dif X^{*}_{1}=[A_{1}X^{*}_{1} +B_{1}(\nabla f_{41})^{-1}(-( B^{\top}_{1}Y_{1}+ \tau\bar{B}_{1}^{\top}\mE[Y_{1}]+D^{\top}_{1}Z_{1}))\\
&\quad\quad\quad\quad\quad\quad\quad\quad+ \bar{B}_{1}\mE[(\nabla f_{41})^{-1}(-( B^{\top}_{1}Y_{1}+ \tau\bar{B}_{1}^{\top}\mE[Y_{1}]+D^{\top}_{1}Z_{1}))]\\
& \quad\quad\quad\quad\quad\quad\quad\quad +\bar{A}_{2}\mE[ X^{*}_{2}]+\bar{A}_{1}\mE[ X^{*}_{1}]+\rho_{1}]\dif t\\
&\quad\quad\quad+\sum^{d}_{j=1}[C_{1j}X^{*}_{1} +D_{1j}(\nabla f_{41})^{-1}(-( B^{\top}_{1}Y_{1}\\
&\quad\quad\quad\quad\quad\quad\quad\quad+ \tau\bar{B}_{1}^{\top}\mE[Y_{1}]+D^{\top}_{1}Z_{1}))+\kappa_{1j}]\dif W_{j}, \\
&\dif Y_{1}=-[\nabla f_{31}(X^{*}_{1})+A_{1}^{\top}Y_{1}
+C_{1}^{\top}Z_{1}\\
&\quad\quad\quad\quad\quad\quad\quad\quad+\bar{A}_{1}^{\top}\mE[Y_{1}]+\bar{A}_{2}^{\top}\mE[Y_{2}] ]\dif t+\sum^{d}_{i=1}Z_{1j}\dif W_{j},\\
&X^{*}_{2}=[A_{2}X^{*}_{2}+B_{2}(\nabla f_{42})^{-1}(-( B^{\top}_{2}Y_{2}+ \tau\bar{B}_{2}^{\top}\mE[Y_{2}]+D^{\top}_{2}Z_{2}))\\
&\quad\quad\quad\quad\quad\quad\quad\quad+ \bar{B}_{2}\mE[(\nabla f_{42})^{-1}(-( B^{\top}_{2}Y_{2}+ \tau\bar{B}_{2}^{\top}\mE[Y_{2}]+D^{\top}_{2}Z_{1}))]\\
&\quad\quad\quad\quad\quad\quad\quad\quad +\bar{A}_{2}\mE[ X^{*}_{1} ]+\rho_{2}]\dif t\\
&\quad\quad\quad+\sum^{d}_{j=1}[C_{2j}X^{*}_{1} +D_{2j}(\nabla f_{41})^{-1}(-( B^{\top}_{2}Y_{2}\\
&\quad\quad\quad\quad\quad\quad\quad\quad+ \tau\bar{B}_{2}^{\top}\mE[Y_{2}]+D^{\top}_{2}Z_{2}))+\kappa_{2j}]\dif W_{j}, \\
&\dif Y_{2}=-[\nabla f_{32}(X^{*}_{2})+A_{2}^{\top}Y_{2}
+C_{2}^{\top}Z_{2}+\bar{A}_{2}^{\top}\mE[Y_{1}] ]\dif t\\
&\quad\quad\quad+\sum^{d}_{j=1}Z_{2j}\dif W_{j}, t\in [0,T],\\
& X^{*}_{1}(0)=\frac{H}{1-\tau^{2}}\bigg\{\bigg.(\nabla f_{11})^{-1}(-\frac{H^{\top}Y_{1}(0)-\tau H^{\top}Y_{2}(0)}{1-\tau^{2}})-\tau\nabla (f_{12})^{-1}(-\frac{H^{\top}Y_{2}(0)-\tau H^{\top}Y_{1}(0)}{1-\tau^{2}})\bigg\}\bigg.+x_{0},\\
 &X^{*}_{2}(0)=\frac{H}{1-\tau^{2}}\bigg\{\bigg.(\nabla f_{12})^{-1}(-\frac{H^{\top}Y_{2}(0)-\tau H^{\top}Y_{1}(0)}{1-\tau^{2}})-\tau\nabla (f_{11})^{-1}(-\frac{H^{\top}Y_{1}(0)-\tau H^{\top}Y_{2}(0)}{1-\tau^{2}})\bigg\}\bigg.+x_{0},\\
& Y_{1}(T)=\nabla f_{21}(X^{*}_{1}(T)+X^{*}_{2}(T))+\nabla f_{22}(X^{*}_{1}(T)+X^{*}_{2}(T)),\\
 & Y_{2}(T)=\nabla f_{21}(X^{*}_{1}(T)+X^{*}_{2}(T))+\nabla f_{22}(X^{*}_{1}(T)+X^{*}_{2}(T)).
\end{cases}
\end{align}
Now, we present our main result as follows:
\bt\label{205}
Let Assumptions 2 and 3 be satisfied. Then, the Hamiltonian system \eqref{37} possesses a unique solution
$$V_{1}(\cdot):=(X_{1}^{*}(\cdot)^{\top},Y_{1}(\cdot)^{\top}, Z_{1}(\cdot)^{\top} )^{\top},  V_{2}(\cdot):=(X^{*}_{2}(\cdot)^{\top},Y_{2}(\cdot)^{\top},Z_{2}(\cdot)^{\top} )^{\top}.$$
Set
$\theta_{i}(t):=((\mE[X^{*}_{1}(t)])^{\top}, (\mE[Y_{1}(t)])^{\top})^{\top}, (\mE[X_{2}(t)])^{\top}, (\mE[Y_{2}(t)])^{\top},   X^{*}_{i}(t)^{\top}, Y_{i}(t)^{\top}, Z_{i}(t)^{\top})^{\top}, i=1,2 . $
Using this solution, the quartet $(\xi^{*}_{1},\xi^{*}_{2},u^{*}_{1},u^{*}_{2})$ defined by \eqref{36}
constitutes the unique optimal controls for Problem (LC)
\et
\begin{proof}
Firstly, we intend to prove the unique solvability of Eq.\eqref{37}.
The Hamiltonian system \eqref{37} takes the form of the System $(\pi)$. To apply Theorem \ref{tt} for establishing the unique solvability of \eqref{37}, it suffices to verify Assumption 1. In what follows, we provide only a detailed verification of the monotonicity condition for the coefficient $\Gamma_{i}(\cdot)$ [cf. Assumption 1(iii)3)], other details are omitted.  Assume that  $\bar{V}_{1}(\cdot):=(\bar{X}_{1}(\cdot)^{\top},\bar{Y}_{1}(\cdot)^{\top}, \bar{Z}_{1}(\cdot)^{\top} )^{\top},  \bar{V}_{2}(\cdot):=(\bar{X}_{2}(\cdot)^{\top},\bar{Y}_{2}(\cdot)^{\top},\bar{Z}_{2}(\cdot)^{\top} )^{\top}$ is a solution of the Hamiltonian system \eqref{37}, $\tilde{V}_{1}(\cdot):=(\tilde{X}_{1}(\cdot)^{\top},\tilde{Y}_{1}(\cdot)^{\top},$ $\tilde{Z}_{1}(\cdot)^{\top} )^{\top},  \tilde{V}_{2}(\cdot):=(\tilde{X}_{2}(\cdot)^{\top},\tilde{Y}_{2}(\cdot)^{\top},\tilde{Z}_{2}(\cdot)^{\top} )^{\top}$ is another solution of the Hamiltonian system \eqref{37} and set
$$\bar{\theta}_{i}(t):=((\mE[\bar{X}_{1}(t)])^{\top}, (\mE[\bar{Y}_{1}(t)])^{\top}), (\mE[\bar{X}_{2}(t)])^{\top}, (\mE[\bar{Y}_{2}(t)])^{\top},  \bar{X}_{i}(t)^{\top}, \bar{Y}_{i}(t)^{\top}, \bar{Z}_{i}(t)^{\top})^{\top}, i=1,2 . $$
$$\tilde{\theta}_{i}(t):=((\mE[\tilde{X}_{1}(t)])^{\top}, (\mE[\tilde{Y}_{1}(t)])^{\top}), (\mE[\tilde{X}_{2}(t)])^{\top}, (\mE[\tilde{Y}_{2}(t)])^{\top},  \tilde{X}_{i}(t)^{\top}, \tilde{Y}_{i}(t)^{\top}, \tilde{Z}_{i}(t)^{\top})^{\top}, i=1,2 . $$
The following result follows from a straightforward calculation, where we have suppressed the dependence on the parameter $t$ and denote $\hat{l}:=\tilde{l}-\bar{l}.$
\begin{align}
&\langle\Gamma_{1}(\tilde{\theta}_{1})-\Gamma_{1}(\bar{\theta}_{1}),\tilde{V}_{1}-\bar{V}_{1}\rangle+\langle\Gamma_{2}(\tilde{\theta}_{2})-\Gamma_{2}(\bar{\theta}_{2}),
\tilde{V}_{2}-\bar{V}_{2}\rangle\no\\
&=-\langle \nabla f_{31}(\tilde{X}_{1})-\nabla f_{31}(\bar{X}_{1}),   \hat{X}_{1}\rangle
-\langle \nabla f_{32}(\tilde{X}_{2})-\nabla f_{32}(\bar{X}_{2}),   \hat{X}_{2}\rangle\no\\
&+\langle\hat{\alpha}_{1}, B^{\top}_{1}\hat{Y}_{1}+\tau\bar{B}^{\top}_{1}\mE[\hat{Y}_{1}]+D^{\top}_{1}\hat{Z}_{1}\rangle+\langle\hat{\alpha}_{2}, B^{\top}_{2}\hat{Y}_{2}+\tau\bar{B}^{\top}_{2}\mE[\hat{Y}_{2}]+D^{\top}_{2}\hat{Z}_{2}\rangle,
\end{align}
where
$$\tilde{\alpha}_{1}:= (\nabla f_{41})^{-1}(-( B^{\top}_{1}\tilde{Y}_{1}+\tau\bar{B}^{\top}_{1}\mE[\tilde{Y}_{1}]+D^{\top}_{1}\tilde{Z}_{1})),$$
$$\tilde{\alpha}_{2}:=(\nabla f_{42})^{-1}(-(B^{\top}_{2}\tilde{Y}_{2}+\tau\bar{B}^{\top}_{2}\mE[\tilde{Y}_{2}]+D^{\top}_{2}\tilde{Z}_{2})), $$
$$ \bar{\alpha}_{1}:= (\nabla f_{41})^{-1}(-(B^{\top}_{1}\bar{Y}_{1}+\tau\bar{B}^{\top}_{1}\mE[\bar{Y}_{1}]+D^{\top}_{1}\bar{Z}_{1})),$$ $$\bar{\alpha}_{2}:=(\nabla f_{42})^{-1}(-(B^{\top}_{2}\bar{Y}_{2}+\tau\bar{B}^{\top}_{2}\mE[\bar{Y}_{2}]+D^{\top}_{2}\bar{Z}_{2})).$$
Given the convexity of $f_{31}, f_{32}$ and uniform convexity of $f_{41}, f_{42},$ it follows from Lemma \ref{100} that
\begin{align*}
-\langle \nabla f_{31}(\tilde{X}_{1})-\nabla f_{31}(\bar{X}_{1}),   \hat{X}_{1}\rangle
-\langle \nabla f_{32}(\tilde{X}_{2})-\nabla f_{32}(\bar{X}_{2}),   \hat{X}_{2}\rangle\leq 0,
\end{align*}
and
\begin{align*}
&\langle  \hat{\alpha}_{1},   B_{1}^{\top}\hat{Y}_{1}+ \tau\bar{B}_{1}^{\top}\mE[\hat{Y}_{1}]+D_{1}^{\top}\hat{Z}_{1}\rangle+
\langle  \hat{\alpha}_{2},   B_{2}^{\top}\hat{Y}_{2}+ \tau\bar{B}_{2}^{\top}\mE[\hat{Y}_{2}]+D_{2}^{\top}\hat{Z}_{2}\rangle\no\\
&\leq -\langle \hat{\alpha}_{1}, \nabla f_{41}(\tilde{\alpha}_{1})- \nabla f_{41}(\bar{\alpha}_{1})   \rangle -\langle \hat{\alpha}_{2}, \nabla f_{42}(\tilde{\alpha}_{2})- \nabla f_{41}(\bar{\alpha}_{2})   \rangle\\
&\leq -\delta|\hat{\alpha}_{1}|^{2} -\delta|\hat{\alpha}_{2}|^{2}.
\end{align*}
We further note that the corresponding monotonicity conditions for $\Gamma_{1}, \Gamma_{2}$ is satisfied by giving the following  definitions:,
\begin{align}
\begin{cases}
&\bar{h}_{11}(v):=\frac{1}{1-\tau^{2}}\bigg.(\nabla f_{11})^{-1}(-v), \bar{h}_{12}(v):=-\frac{\tau}{1-\tau^{2}}\nabla (f_{12})^{-1}(-v),\\
&\bar{h}_{21}(v):=\frac{1}{1-\tau^{2}}\bigg.(\nabla f_{12})^{-1}(-v), \bar{h}_{22}(v):=-\frac{\tau}{1-\tau^{2}}\nabla (f_{11})^{-1}(-v),\\
&  h_{1}(t,u):=(\nabla f_{41})^{-1}(t,-u), h_{2}(t,u):=(\nabla f_{42})^{-1}(t,-u).
\end{cases}
\end{align}

Now, we intend to prove that $(\xi^{*}_{1},\xi^{*}_{2},u^{*}_{1},u^{*}_{2})$  is  the unique quartet of optimal controls. Let $(\xi_{1},\xi_{2},u_{1},u_{2})$
be another arbitrary admissible  quartet. By $(\xi^{*}_{1},\xi^{*}_{2},u^{*}_{1},u^{*}_{2})$ in \eqref{36},   we know that $(\Diamond^{*}_7, \Diamond^{*}_8, \Diamond^{*}_9, \Diamond^{*}_{10}) = (0, 0,0,0).$   Then,   $\eqref{35}$  is reduced to

\begin{align}\label{41}
&J(\xi_{1},\xi_{2}, u_{1}(\cdot), u_{2}(\cdot)) - J(\xi^{*}_{1},\xi^{*}_{2}, u^{*}_{1}(\cdot),u^{*}_{2}(\cdot))\no\\
&\geq \frac{\delta}{2}\bigg\{\bigg. |\xi_{1}+\tau\xi_{2} - \xi^*_{1}-\tau\xi^*_{2}|^2+ |\tau\xi_{1}+\xi_{2} - \tau\xi^*_{1}-\xi^*_{2}|^2\no\\
& +\mE\int^{T}_{0} |u_{1}(t)-u^{*}_{1}(t)|^{2}\dif t +\mE\int^{T}_{0} |u_{2}(t)-u^{*}_{2}(t)|^{2}\dif t    \bigg\}\bigg.> 0.
\end{align}

Thus, $(\xi^{*}_{1},\xi^{*}_{2},u^{*}_{1},u^{*}_{2})$  is the unique quartet of optimal controls for Problem (LC). The proof is thereby established.
\end{proof}

\begin{exa}  Let $n = d = m = k = T = 1$. Consider the following controlled system:
\begin{equation*}
\begin{cases}
&dX_{1}(t) = \frac{1}{1000}\mE[X_{2}(t)]
\dif t+[u_{1}(t) + \sin(W(t))] \dif W(t),\\
&dX_{2}(t) = \frac{1}{1000}\mE[X_{1}(t)] \dif t+[u_{1}(t) + \sin(W(t))] \dif W(t),\\
&X_{1}(0) = \xi_{1},X_{2}(0) = \xi_{2}.
\end{cases}
\end{equation*}
and the following criterion functional:
\begin{equation*}
J(\xi_{1},\xi_{2}, u_{1}(\cdot), u_{2}(\cdot)) = f(\xi_{1})+f(\xi_{2}) + \mathbb{E} \left\{ \frac{1}{4} |X_{1}(1)+X_{2}(1)|^2 + \int_0^1 f(u_{1}(s)) \dif s+ \int_0^1 f(u_{2}(s)) \dif s \right\}.
\end{equation*}
where $f(\cdot)$ is given by \begin{equation}
f_{1}(u) = f_{4}(t,u) = f(u) := \begin{cases} e^{u}-u-1, & u \geq 0 \\ e^{-u}+u-1, & u < 0. \end{cases}
\end{equation}

It is evident that the function $f(\cdot)$ defined above is continuously differentiable.
Specifically,

\begin{equation}\label{201}
\nabla f(u) = \frac{df}{du}(u) = \begin{cases} e^{u}-1, & u \geq 0 \\ 1-e^{-u}, & u < 0. \end{cases}
\end{equation}
Since $(d^2 f)/(du^2)(u) \geq 1$ for any $u \in \mathbb{R},$
we have
$\langle \nabla f(u) - \nabla f(\bar{u}), u - \bar{u} \rangle \geq |u - \bar{u}|^2$ for  any $u,\bar{u}\in \mR.$
Thus, the Lemma \ref{100} in Appendix  shows that $f(\cdot)$ is uniformly convex. Furthermore, since $f(0)=\nabla f(0)=0,$  $Assumption\, 3$ $(3) (4)$ hold.

From \eqref{201}, we have

$$ (\nabla f)^{-1}(u) = \begin{cases}
\ln(1+u), & u \geq 0 \\
-\ln(1-u), & u < 0.
\end{cases} $$
For  $i=1,2,j=1,2,$ set
\begin{align}
\bar{h}_{ij}(u) = h_{i}(t, u) = (\nabla f)^{-1}(-u)
=
\begin{cases}
-\ln(1+u), & u \geq 0, \\
\ln(1-u), & u < 0 ,
\end{cases}
\end{align}
Obviously, we know that
\begin{align}\label{203}
\frac{\partial h_{i}}{\partial u}(t, u) = -\frac{1}{1 + |u|} \in [-1, 0), i=1,2.
\end{align}
On one hand, we assert that the above defined $h_{i}(\cdot,\cdot)$ fulfills $Assumption\, 1(\mathrm{iii}).$ Indeed, \eqref{203} entails the Lipschitz continuity, i.e., the third inequality of $\mathrm{Assumption\, 1\,(iii)\,1)}$ holds true. Moreover, for any $-\infty < \bar{u} < u < \infty$, the Lagrange's mean value theorem implies that there exists a $\tilde{u} \in (\bar{u}, u)$ such that
\begin{equation*}
    h_{i}(t, u) - h_{i}(t, \bar{u}) = \frac{\partial h_{i}}{\partial u}(t, \tilde{u})(u - \bar{u}) \geq -(u - \bar{u})
\end{equation*}
i.e., $u - \bar{u} \geq -(h_{i}(t, u) - h_{i}(t, \bar{u}))$. Therefore,
\begin{equation*}
    (h(t, u) - h(t, \bar{u}))(u - \bar{u}) \leq -|h(t, u) - h(t, \bar{u})|^2.
\end{equation*}
This confirms the fifth inequality of $\mathrm{Assumption\, 1\,(iii)\,1)}$.

On the other hand, also due to \eqref{203} and Lagrange's mean value theorem, there is no constant $c > 0$ such that
\begin{equation*}
    |h(t, u) - h(t, \bar{u})|^2 \geq c|u - \bar{u}|^2 \quad \text{for any } u, \bar{u} \in \mathbb{R}.
\end{equation*}
This indicates that we cannot replace $h$ with $-u$, i.e., the linear version of domination-monotonicity conditions in \cite{Yu2022} does not hold.

Clearly, Assumptions 2 and 3 hold true in this case.

Theorem \ref{205} implies that Problem (LC) with the special setting above admits a unique pair of optimal controls
\begin{align}\label{101}
(\xi^*_{1},\xi^*_{1}, u^*_{1}(\cdot), u^*_{2}(\cdot))= ((\nabla f)^{-1}(-Y_{1}(0)),(\nabla f)^{-1}(-Y_{2}(0))&, (\nabla f)^{-1}(-\frac{1}{1000}\mE[Y_{1}(\cdot)]-(Z_{1}(\cdot)))\no\\
&, (\nabla f)^{-1}(-\frac{1}{1000}\mE[Y_{2}(\cdot)]-(Z_{2}(\cdot)))).
\end{align}
where $(X^*_{i}(\cdot)^{\top}, Y_{i}(\cdot)^{\top}, Z_{i}(\cdot)^{\top})^{\top} \in \mM_F^2(\mathbb{R}^{1+1+1}),i=1,2$ is the unique solution to the following stochastic Hamiltonian system:
\begin{equation}\label{102}
\begin{cases}
&dX^{*}_{1}(t) = \frac{1}{1000}\mE[X^{*}_{2}(t)]\dif t +[(\nabla f)^{-1}(-\frac{1}{1000}\mE[Y_{1}(t)]-Z_{1}(t)) + \sin(W(t))] \dif W(t),\\
&dY_{1}(t) =-\frac{1}{1000}\mE[Y_{2}(t)]\dif t+ Z_{1}(t) \dif W(t),\\
&dX^{*}_{2}(t) = \frac{1}{1000}\mE[X^{*}_{1}(t)]\dif t+[(\nabla f)^{-1}(-\frac{1}{1000}\mE[Y_{1}(t)]-Z_{2}(t))+ \sin(W(t))] \dif W(t),\\
&dY_{2}(t) =-\frac{1}{1000}\mE[Y_{1}(t)]\dif t+ Z_{2}(t) \dif W(t),\\
&X^{*}_{1}(0) = \xi_{1},X^{*}_{2}(0) = \xi_{2}, Y_{1}(1) =\frac{1}{2}( X^{*}_{1}(1)+X^{*}_{2}(1)), Y_{2}(1) =\frac{1}{2}( X^{*}_{1}(1)+X^{*}_{2}(1)).
\end{cases}
\end{equation}

In fact, we can check that the unique solution to \eqref{102} admits the following form:
\begin{equation}\label{103}
\begin{cases}
\xi^{*}_{1}=\xi^{*}_{2},u^{*}_{1}(\cdot)=u^{*}_{2}(\cdot) \\
X^*_{1}(t)=X^*_{2}(t) = Y_{1}(t)= Y_{2}(t) \\
Z_{1}(t)=Z_{2}(t) = u^*_{i}(t) + \sin(W(t)),\\
\mE[X^*_{1}(t)]=\mE[Y_{1}(t)]=\mE[X^*_{2}(t)]=\mE[Y_{2}(t)]=0,
\end{cases}
\end{equation}
for any $(\omega, t) \in \Omega \times [0, 1]$, where the unique optimal pair $(\xi^*_{1},\xi^*_{1}, u^*_{1}(\cdot), u^*_{2}(\cdot))$ satisfies that $\xi_{i}^* = 0,i=1,2,$ and
\begin{equation}
\begin{cases}
e^{u^*_{i}(s)} - 1 + u^*_{i}(s) + \sin(W(s)) = 0, \quad u^*_{i}(s) \geq 0\\
-e^{-u^*_{i}(s)} + 1 + u^*_{i}(s) + \sin(W(s)) = 0, \quad u^*_{i}(s) < 0,
i=1,2.
\end{cases}
\end{equation}

\end{exa}

\section{Application to stochastic linear-quadratic problems with input constraints}
Before giving the main results in this section, we need the following assumptions.\\
$\mathbf{Assumption\, 4:}$
\begin{enumerate}
    \item[(1)] $M_{i} \in \mathbb{S}^m$, $G_{i} \in L_{\mathcal{F}_T}^{\infty}(\mathbb{S}^n)$, $Q_{i}(\cdot) \in L^{\infty}_{\mF}(\mathbb{S}^n)$, and $R_{i}(\cdot) \in L^{\infty}_{\mF}(\mathbb{S}^k), i=1,2$.
    \item[(2)] There exists a constant $\delta > 0$ such that $M - \delta I_m$, $G$, $Q(t)$, and $R(t) - \delta I_k$ are positive semidefinite for almost all $(\omega, t) \in \Omega \times [0, T]$.
\end{enumerate}
In what follows, $I_n$ denotes the $(n \times n)$ identity matrix.

In this section, we turn our attention to the (LQ-IC) problem introduced in Section 1 under Assumptions 2 and 4. Consequently, it follows that
 $$\mJ(\xi_{1},\xi_{2}, u_{1}(\cdot),u_{2}(\cdot))< \infty,\,\mbox{for any}\, (\xi_{1},\xi_{2}, u_{1}(\cdot),u_{2}(\cdot))\in \mR^{m}\times \mR^{m}\times L^{2}_{\mF}(\mR^{k})\times L^{2}_{\mF}(\mR^{k}). $$

Unlike Problem (LC), the controls (inputs) in Problem (LQ-IC) are constrained to closed convex sets. These constraints on $u_i(\cdot), i=1,2$ can be both time-dependent and random. While our approach remains applicable even when the terms involving $X_{i}(T),i=1,2$ and $X_{i}(\cdot),i=1,2$ are replaced by other convex functions possessing suitable properties, we adopt the quadratic criterion functional \eqref{6} in order to present our main idea clearly and concisely. Regarding the constrained sets $U_{0}$ and $U(\cdot)$ (see \cite{DDD8}), the following assumption is introduced.\\
$\mathbf{Assumption\, 5:}$
\begin{enumerate}
    \item[(i)] Both $U_0$ and $U(\omega, t)$ are nonempty, closed, and convex for almost every $(\omega, t) \in \Omega \times [0, T]$.
    \item[(ii)] The indicator function $(\omega, t, u) \mapsto 1_{U(\omega,t)}(u)$ is $\sP \times \sB(\mathbb{R}^k)$-measurable.
    \item[(iii)] There exists a process $a(\cdot) \in L^2_\mathbb{F}(\mathbb{R}^k)$ such that $a(\omega, t) \in U(\omega, t)$ for almost every $(\omega, t) \in \Omega \times [0, T]$.
\end{enumerate}

\br Here, the constraint set $U(\cdot)$ depends on both time and randomness, reflecting the reality that the control constraint limits change according to time and different situations.  We must assume Assumption 5 holds as a compromise. Nonetheless, it is not restrictive and applies to many situations.
Specifically, Assumption 5 obviously holds when $U$ does not depend on $(\omega, t)$ and Assumption 5(i) holds.
Furthermore, when $U(\omega, t)$ is uniformly bounded for almost all $(\omega, t) \in \Omega \times [0, T]$, Assumption 5(iii) is superfluous.  In fact, under Assumption 5(ii),

$$
\Lambda = \{ (\omega, t, u) \in \Omega \times [0, T] \times \mathbb{R}^k | 1_{U(\omega, t)}(u) = 1 \}
$$
is $\sP \otimes \mathcal{B}(\mathbb{R}^k)$-measurable.  Consequently, according to the measurable selection theorem, there is a $\mathcal{P}$-measurable process $a(\cdot)$ such that $a(\omega, t) \in U(\omega, t)$.
\er
Assumption 5(iii) guarantees that the set $\mathcal{U}$ defined by \eqref{as} is a non-empty set.  We refer to $U_0 \times \mathcal{U}$ as the admissible control set.  If $(\xi_{1}, \xi_{2}, u_{1}(\cdot), u_{2}(\cdot)) \in U_0 \times U_0 \times \mathcal{U}\times \mathcal{U}$, then it is termed a pair of admissible controls.  In this situation, $X_{1}(\cdot) \equiv X_{1}(\cdot ; \xi_{1}, u_{1}(\cdot)),X_{2}(\cdot) \equiv X_{2}(\cdot ; \xi_{2}, u_{2}(\cdot))$ and $(\xi_{1},\xi_{1},u_{1}(\cdot), u_{2}(\cdot), x_{1}(\cdot), x_{2}(\cdot))$ are referred to as the corresponding admissible state and an admissible sextet, respectively.
To define the optimal control pair for Problem (LQ-IC), we introduce the Hamiltonian system below, keeping in mind that the argument $t$ is omitted for brevity:

\begin{equation}\label{43}
\begin{cases}
&\dif X^*_{1} = [A_{1}X^*_{1} + B_{1}\Pi(-R_{1}^{-1}( B^{\top}_{1}Y_{1}+ \tau\bar{B}_{1}^{\top}\mE[Y_{1}]+D^{\top}_{1}Z_{1})) \\
&\quad \quad \quad\quad \quad \quad+\bar{B}_{1}\mE[\Pi(-R_{1}^{-1}( B^{\top}_{1}Y_{1}+ \tau\bar{B}_{1}^{\top}\mE[Y_{1}]+D^{\top}_{1}Z_{1})) ]\\
&\quad \quad \quad\quad \quad \quad+\bar{A}_{2}\mE[ X^{*}_{2}]+\bar{A}_{1}\mE[ X^{*}_{1}]+\rho_{1}]\dif t\\
 &\quad \quad \quad+ \sum_{i=1}^{d}[C_{1j}X^*_{1} + D_{1j}\Pi(-R^{-1}_{1}( B^{\top}_{1}Y_{1}+ \tau\bar{B}_{1}^{\top}(t)\mE[Y_{1}]+D^{\top}_{1}Z_{1}))+\kappa_{1j}  ]\dif W_j \\
&\dif Y_{1} = -[Q_{1}X^*_{1} + A^{\top}_{1}Y_{1} + C^{\top}Z_{1}+\bar{A}^{\top}_{2}\mE[ Y_{2}]+\bar{A}^{\top}_{1}\mE[ Y_{1}]]dt + \sum_{i=j}^{d} Z_{1j} \dif W_j \\
&\dif X^*_{2} = [A_{2}X^*_{2} + B_{2}\Pi(-R_{2}^{-1}( B^{\top}_{2}Y_{2}+ \tau\bar{B}_{2}^{\top}\mE[Y_{2}]+D^{\top}_{2}Z_{2}))\\
&\quad \quad \quad\quad \quad \quad+\bar{B}_{2}\mE[\Pi(-R_{2}^{-1}( B^{\top}_{2}Y_{2}+ \tau\bar{B}_{2}^{\top}\mE[Y_{2}]+D^{\top}_{2}Z_{2})) ]+\bar{A}_{2}\mE[ X^{*}_{1}]+\rho_{2}]\dif t\\
 &\quad \quad \quad+ \sum_{i=1}^{d}[C_{2j}X^*_{1} + D_{2j}\Pi(-R^{-1}_{2}( B^{\top}_{2}Y_{2}+ \tau\bar{B}_{2}^{\top}\mE[Y_{2}]+D^{\top}_{2}Z_{2}))+\kappa_{2 j}  ]\dif W_j \\
&\dif Y_{2} = -[Q_{2}X^*_{2} + A_{2}^{\top}Y_{2} + C^{\top}_{2}Z_{2}+\bar{A}_{2}\mE[ Y_{1}]]dt + \sum_{i=1}^{d} Z_{2i} \dif W_i \\
&X^*_{1}(0) = H\Pi_0(-M^{-1}_{1}H^{\top}Y_{1}(0)) + x_0, \\
&  Y_{1}(T) = G_{1}(X^*_{1}(T)+X^*_{2}(T))+G_{2}(X^*_{1}(T)+X^*_{2}(T)),\\
&X^*_{2}(0) = H\Pi_0(-M^{-1}_{2}H^{\top}Y_{2}(0)) + x_0, \\
& Y_{2}(T) = G_{1}(X^*_{1}(T)+X^*_{2}(T))+G_{2}(X^*_{1}(T)+X^*_{2}(T)).
\end{cases}
\end{equation}

Here, $\Pi_0(\cdot)$ and $\Pi(\omega, t, \cdot)$ represent the projections from $\mathbb{R}^m$ and $\mathbb{R}^k$, respectively, onto their corresponding closed convex subsets $U_0$ and $U(\omega, t)$, under the specified norm

$$
|\cdot|_{M_{i}}:= \langle \cdot, \cdot \rangle_{M_{i}}^{1/2} := \langle M_{i} \cdot, \cdot \rangle^{1/2}, i=1,2,
$$
and (respectively)
$$
|\cdot|_{R_{i}(\omega, t)} := \langle \cdot, \cdot \rangle_{R_{i}(\omega, t)}^{1/2} := \langle R_{i}(\omega, t) \cdot, \cdot \rangle^{1/2}, i=1,2
$$
for almost every $(\omega, t) \in \Omega \times [0, T]$.
Before giving our main result in this section,  we present the following lemma, which can be found in Lemma 3 of  \cite{LiuNiuWangYu2026}.
  \bl\label{d}
  Under Assumptions 4 and 5, both $\Pi_0(\cdot)$ and $\Pi(t,\cdot)$ are Lipschitz continuous uniformly with respect to almost all $(\omega, t) \in \Omega \times [0, T]$. Moreover, $\Pi$ is $\mathcal{P} \otimes \mathcal{B}(\mathbb{R}^k)$-measurable and the following inequality holds
  \begin{align}\label{e}
  \delta |u|^2 \leq |u|_{R_{i}(\omega,t)}^2 \leq \|R_{i}(\cdot)\|_{L^{\infty}_{\mF}(\mS^{k})} |u|^2 \text{ for any } u \in \mathbb{R}^k.
\end{align}
  \begin{align}\label{e2}
  |\Pi(\omega,t,u)-\Pi(\omega,t,\bar{u})|^{2}\leq \frac{\|R_{i}(\cdot)\|_{L^{\infty}_{\mF}(\mS^{k})}}{\delta}|u-\bar{u}|^{2}.
\end{align}
 \el

Now, we give the main result in this section.

 \bt
 Suppose Assumptions 2, 4, and 5 are satisfied. Then, the Hamiltonian system \eqref{43} possesses a unique solution $ (V_{1}(\cdot),V_{2}(\cdot))$  where  $V_{1}(\cdot) := (X^*_{1}(\cdot)^{\top}, Y_{1}(\cdot)^{\top}, Z_{1}(\cdot)^{\top})^{\top} \in \mM_{\mathbb{F}}^2(\mathbb{R}^{n+n+nd})$, $V_{2}(\cdot) := (X^*_{2}(\cdot)^{\top}, Y_{2}(\cdot)^{\top}, Z_{2}(\cdot)^{\top})^{\top} \in \mM_{\mathbb{F}}^2(\mathbb{R}^{n+n+nd})$.  Moreover,

\begin{equation}\label{44+-+}
\begin{cases}
\xi^*_{1} =  \Pi_0(-M^{-1}_{1}H^{\top}Y_{1}(0)), \\
\xi^*_{2} = \Pi_0(-M^{-1}_{2}H^{\top}Y_{2}(0)), \\
u^*_{1}(\cdot) = \Pi(\cdot, -R_{1}(\cdot)^{-1}(B^{\top}_{1}(\cdot)Y_{1}(\cdot)+ \tau\bar{B}_{1}^{\top}(\cdot)\mE[Y_{1}(\cdot)]+D^{\top}_{1}(\cdot)Z_{1}(\cdot) ) ),\\
u^*_{2}(\cdot)  = \Pi(\cdot, -R_{2}(\cdot)^{-1}(B^{\top}_{2}(\cdot)Y_{2}(\cdot)+ \tau\bar{B}_{2}^{\top}(\cdot)\mE[Y_{2}(\cdot)]+D^{\top}_{2}(\cdot)Z_{2}(\cdot) )
\end{cases}
\end{equation}
is the quartet  of optimal controls for Problem (LQ-IC).
\et
\begin{proof}
$\mathbf{Stage\, 1}:$ we intend to prove  the unique solvability of \eqref{43}.

To utilize Theorem \ref{tt} to demonstrate the unique solvability of \eqref{43}, it is necessary to confirm Assumption 1. We present exclusively the detailed confirmations of $b_{i}(\cdot, 0) \in L_{\mF}^{2}(\mathbb{R}^n), i=1,2$ [refer to Assumption 1(i)] and the monotonicity property of $\Gamma_{i}(\cdot),i=1,2$ [refer to Assumption 1(iii)3)]. Further specifics are not shown.

Initially, in this specific scenario, $b_{1}(\cdot, 0) = B_{1}(\cdot)\Pi(\cdot, 0) + \rho_{1}(\cdot), b_{2}(\cdot, 0) = B_{2}(\cdot)\Pi(\cdot, 0) + \rho_{2}(\cdot)$. Given that $B_{1}(\cdot), B_{2}(\cdot) \in L^{\infty}(\mathbb{R}^{n \times k})$ and $\rho_{1}(\cdot), \rho_{2}(\cdot) \in L_{\mF}^{2}(\mathbb{R}^n)$, the task simplifies to verifying $\Pi(\cdot, 0) \in L_{\mF}^{2}(\mathbb{R}^k)$.  Indeed, with the assistance of Assumption 5(iii) and the Lipschitz continuity of $\Pi(t, \cdot)$ (as stated in Lemma \ref{d}), we obtain

\begin{align*}
|\Pi(t, 0)| &\leq |\Pi(t, 0) - a(t)| + |a(t)| \\
&= |\Pi(t, 0) - \Pi(t, a(t))| + |a(t)| \\
&\leq \sqrt{\frac{\|R(\cdot)\|_{L^{\infty}_{\mF}(\mS^{k})}}{\delta}} |a(t)| + |a(t)|
\end{align*}
for almost every $(\omega, t) \in \Omega \times [0, T]$. This leads to the conclusion that $\Pi(\cdot, 0) \in L_{\mF}^2(\mathbb{R}^k) $.

Secondly, assume that $(\tilde{V}_{1}(\cdot),\tilde{V}_{2}(\cdot))$ and $(\bar{V}_{1}(\cdot),\bar{V}_{2}(\cdot))$ are two solution of  Eq.\eqref{43}  ,where
 $$\tilde{V}_{1}(\cdot) := (\tilde{X}_{1}(\cdot)^{\top}, \tilde{Y}_{1}(\cdot)^{\top}, \tilde{Z}_{1}(\cdot)^{\top})^{\top} \in \mM_{\mathbb{F}}^2(\mathbb{R}^{n+n+nd}),$$  $$\tilde{V}_{2}(\cdot) := (\tilde{X}_{2}(\cdot)^{\top}, \tilde{Y}_{2}(\cdot)^{\top}, \tilde{Z}_{2}(\cdot)^{\top})^{\top} \in \mM_{\mathbb{F}}^2(\mathbb{R}^{n+n+nd})$$ and $$\bar{V}_{1}(\cdot) := (\bar{X}_{1}(\cdot)^{\top}, \bar{Y}_{1}(\cdot)^{\top}, \bar{Z}_{1}(\cdot)^{\top})^{\top} \in \mM_{\mathbb{F}}^2(\mathbb{R}^{n+n+nd}),$$ $$\bar{V}_{2}(\cdot) := (\bar{X}_{2}(\cdot)^{\top}, \bar{Y}_{2}(\cdot)^{\top}, \bar{Z}_{2}(\cdot)^{\top})^{\top} \in \mM_{\mathbb{F}}^2(\mathbb{R}^{n+n+nd}).$$    A direct computation yields (where the argument $t$ is omitted and $\hat{l} := \tilde{l} - \bar{l}$)

\begin{align*}
&\langle \Gamma_{1}(\tilde{\theta}_{1}) - \Gamma_{1}(\bar{\theta}_{1}), \hat{V_{1}} \rangle +\langle \Gamma_{2}(\tilde{\theta}_{2}) - \Gamma_{2}(\bar{\theta}_{2}), \hat{V_{2}} \rangle \\
&= - \langle Q_{1} \hat{X}_{1}, \hat{X}_{1} \rangle + \langle \hat{\varrho}_{1},  B^{\top}_{1}\hat{Y}_{1}+ \tau\bar{B}_{1}^{\top}(t)\mE[\hat{Y}_{1}]+D^{\top}_{1}\hat{Z}_{1} \rangle \\
&- \langle Q_{2} \hat{X}_{2}, \hat{X}_{2} \rangle + \langle \hat{\varrho}_{1},  B^{\top}_{2}\hat{Y}_{1}+ \tau\bar{B}_{2}^{\top}(t)\mE[\hat{Y}_{2}]+D^{\top}_{2}\hat{Z}_{2} \rangle\\
&\leq \langle \hat{\varrho}_{1}, B^{\top}_{1}\hat{Y}_{1}+ \tau\bar{B}_{1}^{\top}\mE[\hat{Y}_{1}]+D^{\top}_{1}\hat{Z}_{1}  \rangle\no\\
&+ \langle \hat{\varrho}_{2}, B^{\top}_{2}\hat{Y}_{1}+ \tau\bar{B}_{2}^{\top}\mE[\hat{Y}_{2}]+D^{\top}_{2}\hat{Z}_{2} \rangle,
\end{align*}
where \begin{equation}\label{44}
\begin{cases}
\tilde{\varrho}_{1} = \Pi( -R_{1}^{-1}(B^{\top}_{1}\tilde{Y}_{1}+ \tau\bar{B}_{1}^{\top}\mE[\tilde{Y}_{1}]+D^{\top}_{1}\tilde{Z}_{1} ) )\\
\bar{\varrho}_{1}  = \Pi( -R_{1}^{-1}(B^{\top}_{1}\bar{Y}_{1}+ \tau\bar{B}_{1}^{\top}\mE[\bar{Y}_{1}]+D^{\top}_{1}\bar{Z}_{1} )\\
\tilde{\varrho}_{2}  = \Pi( -R_{2}^{-1}(B^{\top}_{2}\tilde{Y}_{2}+ \tau\bar{B}_{2}^{\top}\mE[\tilde{Y}_{2}]+D^{\top}_{2}\tilde{Z}_{2} )\\
\bar{\varrho}_{2} = \Pi( -R_{2}^{-1}(B^{\top}_{2}\bar{Y}_{2}+ \tau\bar{B}_{2}^{\top}\mE[\bar{Y}_{2}]+D^{\top}_{2}\bar{Z}_{2} ).
\end{cases}
\end{equation}
By  Lemma \ref{c} and \eqref{e} in Lemma \ref{d},
\begin{align*}
    &\langle \hat{\varrho}_{1}, B^{\top}_{1}\hat{Y}_{1}+ \tau\bar{B}_{1}^{\top}\mE[\hat{Y}_{1}]+D^{\top}_{1}(t)\hat{Z}_{1} \rangle\\
    & = - \langle \hat{\varrho}_{1}, -R_{1}^{-1} ( B^{\top}_{1}\hat{Y}_{1}+ \tau\bar{B}_{1}^{\top}(t)\mE[\hat{Y}_{1}]+D^{\top}_{1}(t)\hat{Z}_{1} ) \rangle_{R_{1}} \leq -|\hat{\varrho}_{1}|^2_{R_{1}} \leq -\delta|\hat{\varrho}_{1}|^2.
\end{align*}
and
\begin{align*}
    &\langle \hat{\varrho}_{2},B^{\top}_{2}\hat{Y}_{2}+ \tau\bar{B}_{2}^{\top}\mE[\hat{Y}_{2}]+D^{\top}_{2}(t)\hat{Z}_{2}\rangle\\
    & = - \langle \hat{\varrho}_{2}, -R_{2}^{-1} (B^{\top}_{2}\hat{Y}_{2}+ \tau\bar{B}_{2}^{\top}(t)\mE[\hat{Y}_{2}]+D^{\top}_{2}(t)\hat{Z}_{2}) \rangle_{R_{2}} \leq -|\hat{\varrho}_{2}|^2_{R_{2}} \leq -\delta|\hat{\varrho}_{2}|^2.
\end{align*}
Based on the following definitions,

\begin{equation*}
\begin{cases}
    \bar{h}_{11}(v) =   \Pi_0 (-M_{1}^{-1}v),\\
    \bar{h}_{12}(v) =   0,\\
    \bar{h}_{21}(v) =  \Pi_0 (-M_{2}^{-1}v),\\
    \bar{h}_{22}(v) =  0,\\
    h_{1}(t, u) = \Pi (t, -R_{1}(t)^{-1}u)\\
  h_{2}(t, u) = \Pi (t, -R_{2}(t)^{-1}u).
\end{cases}
\end{equation*}
we can infer that the associated monotonicity property is valid for $\Gamma_{1}(\cdot), \Gamma_{2}(\cdot)$.

$\mathbf{Stage\, 2}:$ $(\xi^*_{1},\xi^*_{2}, u^*_{1}(\cdot), u^*_{2}(\cdot))$ is the only optimal control quartet.
Assume that the corresponding admissible sextet is $(\xi^*_{1},\xi^*_{2}, u^*_{1}(\cdot), u^*_{2}(\cdot),X^{*}_{1}(\cdot), X^{*}_{2}(\cdot)).$\\
Suppose $(\xi_{1},\xi_{1}, u_{1}(\cdot), u_{2}(\cdot), X_{1}(\cdot),X_{2}(\cdot))$ is another arbitrary admissible sextet.  We examine the difference

\begin{equation*}
    \mathcal{J}(\xi_{1},\xi_{1}, u_{1}(\cdot), u_{2}(\cdot)) - \mathcal{J}(\xi^*_{1},\xi^*_{1}, u^*_{1}(\cdot), u^*_{2}(\cdot)) = \Re_{81} + \Re_{82}+\Re_{91} + \Re_{92}
\end{equation*}
where (the time argument $t$ is omitted for brevity)

\begin{align*}
    \Re_{81} &:= \frac{1}{2}\langle M_{1}(\xi_{1} - \xi^{*}_{1}), \xi_{1} - \xi^{*}_{1} \rangle \\
    &+\frac{1}{2} \mathbb{E}\bigg\{\bigg. \langle G_{1}(X_{1}(T)+X_{2}(T) - X^*_{1}(T))- X^*_{2}(T)),X_{1}(T)+X_{2}(T) - X^*_{1}(T))- X^*_{2}(T) \rangle\\ + &\int_{0}^{T} [ \langle Q_{1}(X_{1} - X^{*}_{1}), X_{1} - X^{*}_{1} \rangle + \langle R_{1}(u_{1} - u^{*}_{1}), u_{1} - u^{*}_{1} \rangle ] \dif t\bigg\}\bigg..
\end{align*}

\begin{align*}
\Re_{91} &:= \langle M_{1}\xi^{*}_{1}, \xi_{1}-\xi^{*}_{1} \rangle \\
&+  \mathbb{E}\bigg\{\bigg. \langle G_{1}( X^*_{1}(T))+X^*_{2}(T)),X_{1}(T)+X_{2}(T) - X^*_{1}(T))- X^*_{2}(T) \rangle\\
&\quad + \int_0^T \left[ \langle Q_{1}X^*_{1}, X_{1} - X^*_{1} \rangle + \langle R_{1}u^*_{1}, u_{1} - u^*_{1} \rangle \right] dt \bigg\}.
\end{align*}

\begin{align*}
    \Re_{82} &:= \frac{1}{2}\langle M_{2}(\xi_{2}- \xi^{*}_{2}), \xi_{2}- \xi^{*}_{2} \rangle\\\\
    &+ \frac{1}{2} \mathbb{E}\bigg\{\bigg. \langle G_{2}(X_{1}(T)+X_{2}(T) - X^*_{1}(T))- X^*_{2}(T)),X_{1}(T)+X_{2}(T) - X^*_{1}(T))- X^*_{2}(T) \rangle\\ + &\int_{0}^{T} [ \langle Q_{2}(X_{2} - X^{*}_{2}), X_{2} - X^{*}_{2} \rangle + \langle R_{2}(u_{2} - u^{*}_{2}), u_{2} - u^{*}_{2} \rangle ] \dif t\bigg\}\bigg..
\end{align*}

\begin{align*}
\Re_{92} &:= \langle M_{2}\xi^{*}_{2}, \xi_{2} - \xi^{*}_{2} \rangle  \\
& + \mathbb{E}\bigg\{\bigg. \langle G_{2}( X^*_{1}(T))+ X^*_{2}(T)),X_{1}(T)+X_{2}(T) - X^*_{1}(T))- X^*_{2}(T) \rangle\\
&\quad + \int_0^T \left[ \langle Q_{2}X^*_{2}, X_{2} - X^*_{2} \rangle + \langle R_{2}u^*_{2}, u_{2} - u^*_{2} \rangle \right] \dif t \bigg\}.
\end{align*}
Initially, applying It\^{o}'s formula to $\langle Y_{1}(\cdot), X_{1}(\cdot) - X^*_{1}(\cdot) \rangle$ and $\langle Y_{2}(\cdot), X_{2}(\cdot) - X^*_{2}(\cdot) \rangle$ results in
\begin{align*}
&\mathbb{E}\left\{ \langle G_{1}X^*_{1}(T), X_{1}(T) - X^*_{1}(T) \rangle + \int_0^T \langle Q_{1}X^*_{1}, X_{1} - X^*_{1} \rangle \dif t  \right\}\\
&+\mathbb{E}\left\{ \langle G_{2}X^*_{2}(T), X_{2}(T) - X^*_{2}(T) \rangle + \int_0^T \langle Q_{2}X^*_{2}, X_{2} - X^*_{2} \rangle \dif t\right\}\\
&= \langle H^T Y_{1}(0), \xi_{1} - \xi^*_{1} \rangle +\langle H^T Y_{2}(0), \xi_{2} - \xi^*_{2} \rangle\\
 &\quad +\mE\int^{T}_{0}\langle   B^{\top}_{1}Y_{1}+ \tau\bar{B}_{1}^{\top}\mE[Y_{1}]+D^{\top}_{1}Z_{1}, u_{1}-u^{*}_{1}   \rangle\dif t\no\\
&\quad +\mE\int^{T}_{0}\langle B^{\top}_{2}Y_{2}+ \tau\bar{B}_{2}^{\top}\mE[Y_{2}]+D^{\top}_{2}Z_{2},   u_{2}-u^{*}_{2} \rangle\dif t.
\end{align*}
Substituting this into $\Re_{91}+\Re_{92}$ gives
\begin{align*}
&\Re_{91}+\Re_{92} = \langle M_{1}\xi^{*}_{1}+H^T Y_{1}(0), \xi_{1}-\xi^{*}_{1} \rangle+\langle M_{2}\xi^{*}_{2}+H^T Y_{2}(0), \xi_{2} - \xi^{*}_{2} \rangle \\
& + \mathbb{E} \int_0^T \langle R_{1}u^*_{1} + B^{\top}_{1}Y_{1}+ \tau\bar{B}_{1}^{\top}\mE[Y_{1}]+D^{\top}_{1}Z_{1}, u_{1} - u^*_{1} \rangle dt\no\\
&+ \mathbb{E} \int_0^T \langle R_{2}u^*_{2} + \tau\bar{B}_{2}^{\top}\mE[Y_{2}]+D^{\top}_{2}Z_{2}, u_{2} - u^*_{2} \rangle \dif t.
\end{align*}
Because of \eqref{a} in Lemma \ref{c} and the definition of $\xi^*$ [refer to \eqref{44}], The sum of the first two terms of $\Re_{91}+\Re_{92}$  simplifies to  \begin{align*}
&\langle M_{1}\xi^{*}_{1}+H^T Y_{1}(0), \xi_{1}-\xi^{*}_{1} \rangle+\langle M_{2}\xi^{*}_{2}+H^T Y_{2}(0), \xi_{2} - \xi^{*}_{2} \rangle\\
&= - \langle -M^{-1}_{1} H^T Y_{1}(0) - \xi^*_{1}, \xi_{1} - \xi^*_{1} \rangle_{M_{1} }- \langle -M^{-1}_{2} H^T Y_{2}(0) - \xi^*_{2}, \xi_{2} - \xi^*_{2} \rangle_{M_{2} }\geq 0.
\end{align*}
In the same way, the second term of $\Re_{91}+\Re_{92}$ is also greater than or equal to zero. As a consequence, we obtain $\Re_{91}+\Re_{92} \geq 0$.

On the other hand, by $\mathrm{Assumption\, 4\,(2})$ and $(\xi_{1},\xi_{1}, u_{1}(\cdot), u_{2}(\cdot))\neq (\xi^*_{1},\xi^*_{1}, u^*_{1}(\cdot), u^*_{2}(\cdot)),$ it yields

\begin{align*}
&\mathcal{J}(\xi, u(\cdot)) - \mathcal{J}(\xi^{*}, u^{*}(\cdot)) \geq \Re_{81}+\Re_{82}\\
&\geq \frac{\delta}{2} \left\{ |\xi_{1}- \xi^{*}_{1}|^{2} + \mathbb{E} \int_{0}^{T} |u_{1} - u^{*}_{1}|^{2} dt \right\} \\
&+\frac{\delta}{2} \left\{ |\xi_{2} - \xi^{*}_{2}|^{2} + \mathbb{E} \int_{0}^{T} |u_{2} - u^{*}_{2}|^{2} dt \right\} > 0
\end{align*}
 Consequently, $(\xi^*_{1},\xi^*_{1}, u^*_{1}(\cdot), u^*_{2}(\cdot))$ defined by \eqref{44+-+} constitutes the unique optimal control  quadruple for Problem (LQ-IC). We complete the proof.

\end{proof}

\section{Appendix}
\subsection{The proof of the results  in Sec. 3.}

\begin{proof}[The proof of Lemma \ref{t1}]
  First,  we intend to  prove \eqref{11}.  Employing a standard estimate for stochastic differential equations (see, e.g., \cite{Zhang2017}), we obtain
\begin{align}\label{13}
&\mE[\sup_{0\leq t\leq T}|X_{1}(t)|^{2}] + \mE[\sup_{0\leq t\leq T}|X_{2}(t)|^{2}] \no\\
&\leq C\bigg\{\bigg.   |\Psi_{1}(Y_{1}(0),Y_{2}(0))-\Psi_{1}(0,Y_{2}(t))|^{2}
+|\Psi_{1}(0,Y_{2}(0))-\Psi_{1}(0,0)|^{2}\no\\
&\quad \quad \quad\quad \quad +|\Psi_{2}(Y_{1}(0),Y_{2}(0))-\Psi_{2}(Y_{1}(t), 0)|^{2}
+|\Psi_{2}(Y_{1}(0),0)-\Psi_{2}(0,0)|^{2}\no\\
&\quad \quad \quad\quad \quad +|\Psi_{1}(0,0)|^{2}+ |\Psi_{2}(0,0)|^{2}\no\\
&   +\mE\bigg[\bigg.\bigg(\bigg.\int^{T}_{0}|b_{1}(t,0)|\dif t  \bigg)\bigg.^{2}+\bigg(\bigg.\int^{T}_{0}|b_{2}(t,0)|\dif t  \bigg)\bigg.^{2}  +\bigg.\bigg(\bigg.\int^{T}_{0}|\sigma_{1}(t, 0)|\dif t  \bigg)\bigg.^{2}+\bigg(\bigg.\int^{T}_{0}|\sigma_{2}(t,0)|\dif t  \bigg)\bigg.^{2} \no \\
&   +\mE\bigg[\bigg.\bigg(\bigg.\int^{T}_{0}|b_{1}(t,0, \mE[Y_{2}(t)],0)-b_{1}(t,0)|\dif t  \bigg)\bigg.^{2} +\bigg(\bigg.\int^{T}_{0}|b_{2}(t,0, \mE[Y_{1}(t)], 0)-b_{2}(t,0)|\dif t  \bigg)\bigg.^{2}  \no\\
&    +\bigg.\bigg(\bigg.\int^{T}_{0}|\sigma_{1}(t,0, \mE[Y_{2}(t)], 0)-\sigma_{1}(t, 0)|\dif t  \bigg)\bigg.^{2} +\bigg(\bigg.\int^{T}_{0}|\sigma_{2}(t,0,  \mE[Y_{1}(t)], 0)-\sigma_{2}(t, 0)|\dif t  \bigg)\bigg.^{2} \no \\
 &   +\bigg(\bigg.\int^{T}_{0}|b_{1}(t, \Upsilon_{1}(t) )-\mE[b_{1}(t,0, \mE[Y_{2}(t)], 0)]|\dif t  \bigg)\bigg.^{2}\no\\
 &    +\bigg(\bigg.\int^{T}_{0}|b_{2}(t,\Upsilon_{2}(t)) -\mE[b_{2}(t,0,  \mE[Y_{1}(t)], 0)]|\dif t  \bigg)\bigg.^{2}   \no\\
&  +\int^{T}_{0}|\mE[\sigma_{1}(t,\Upsilon_{1}(t))-\sigma_{1}(t,0, \mE[Y_{2}(t)], 0)]|^{2}\dif t \no\\
 &   +\int^{T}_{0}|\sigma_{2}(t,\Upsilon_{2}(t))-\sigma_{2}(t,0, \mE[Y_{1}(t)], 0)]|^{2}\dif t  \bigg]\bigg.  \bigg\}\bigg.,
\end{align}
where $\Upsilon_{i}(t):=(0^{\top}, (\mE[Y_{1}(t)])^{\top}, 0^{\top}, (\mE[Y_{2}(t)])^{\top},   0^{\top}, Y_{i}(t)^{\top}, Z_{i}(t)^{\top})^{\top}, i=1,2. $
Then, by domination conditions and Lipschitz  conditions in Assumption 1 (also noting that $\epsilon>0, \varepsilon>0$ are small enough),  we have
\begin{align}\label{14}
&\mE[\sup_{0\leq t\leq T}|X_{1}(t)|^{2}] + \mE[\sup_{0\leq t\leq T}|X_{2}(t)|^{2}]\no \\
&\leq K\bigg\{\bigg.I_{1}+ |\Psi_{1}(0,0)|^{2}+ |\Psi_{2}(0,0)|^{2}+ (\varepsilon+\tau) \mE[Y_{1}(t)]+ (\varepsilon+\tau) \mE[Y_{2}(t)]\no\\
&\quad \quad   +\mE\bigg[\bigg.\bigg(\bigg.\int^{T}_{0}|b_{1}(t,0)|\dif t  \bigg)\bigg.^{2} +\bigg(\bigg.\int^{T}_{0}|b_{2}(t, 0)|\dif t  \bigg)\bigg.^{2}  \bigg]\bigg.\no \\
& \quad \quad   +\mE\bigg[\bigg.\bigg(\bigg.\int^{T}_{0}|\sigma_{1}(t,0)|\dif t  \bigg)\bigg.^{2} +\bigg(\bigg.\int^{T}_{0}|\sigma_{2}(t, 0)|\dif t  \bigg)\bigg.^{2}  \bigg]\bigg. \bigg\}\bigg.,\no \\
\end{align}
where \begin{align}\label{15}
I_{1}&:=\bigg|\bigg.\bar{h}_{11}(\frac{H^{\top}Y_{1}(0)-\tau H^{\top}Y_{2}(0)}{1-\tau^{2}})-\bar{h}_{11}(\frac{-\tau H^{\top}Y_{2}(0)}{1-\tau^{2}})\bigg|\bigg.^{2}\no\\
&\quad+\bigg|\bigg.\bar{h}_{12}(\frac{H^{\top}Y_{2}(0)-\tau H^{\top}Y_{1}(0)}{1-\tau^{2}})-\bar{h}_{12}(\frac{ H^{\top}Y_{2}(0)}{1-\tau^{2}})\bigg|\bigg.^{2}\no\\
&\quad+\bigg|\bigg.\bar{h}_{21}(\frac{H^{\top}Y_{2}(0)-\tau H^{\top}Y_{1}(0)}{1-\tau^{2}})-\bar{h}_{21}(\frac{ -\tau H^{\top}Y_{1}(0)}{1-\tau^{2}})\bigg|\bigg.^{2}\no\\
&\quad+\bigg|\bigg.\bar{h}_{22}(\frac{H^{\top}Y_{1}(0)-\tau H^{\top}Y_{2}(0)}{1-\tau^{2}})-\bar{h}_{22}(\frac{ H^{\top}Y_{1}(0)}{1-\tau^{2}})\bigg|\bigg.^{2}\no\\
&+\mE\int^{T}_{0}|h_{1}(t,B_{1}(t)^{\top}Y_{1}(t)+\tau\bar{B}_{1}(t)^{\top}\mE[Y_{1}(t)]+D_{1}(t)^{\top}Z_{1}(t))-h_{1}(t,0)|^{2}\dif t\no\\
&+\mE\int^{T}_{0}|h_{2}(t,B_{2}(t)^{\top}Y_{2}(t)+\tau\bar{B}_{2}(t)^{\top}\mE[Y_{2}(t)]+D_{2}(t)^{\top}Z_{2}(t))-h_{2}(t,0)|^{2}\dif t.
\end{align}
Furthermore, by the basic estimate of BSDEs(see, e.g., \cite{PardouxPeng1990}) and the Lipschitz conditions,  it holds that
\begin{align}\label{16}
&\mE\bigg[\bigg. \sup_{0\leq t\leq T}|Y_{1}(t)|^{2}+\sup_{0\leq t\leq T}|Y_{2}(t)|^{2}+ \int^{T}_{0}|Z_{1}(t)|^{2}\dif t+\int^{T}_{0}|Z_{2}(t)|^{2}\dif t\bigg]\bigg. \no\\
&\leq C\mE\bigg[\bigg. |\Phi_{1}(0,0)|^{2}+ |\Phi_{2}(0,0)|^{2}
   +\bigg[\bigg.\bigg(\bigg.\int^{T}_{0}|f_{1}(t,0)|\dif t  \bigg)\bigg.^{2} \no\\
   & \quad \quad \quad\quad \quad \quad+\bigg(\bigg.\int^{T}_{0}|f_{2}(t, 0)|^{2}\dif t  \bigg)\bigg.^{2}+   \sup_{0\leq t\leq T}|X_{1}(t)|^{2}+\sup_{0\leq t\leq T}|X_{2}(t)|^{2}\bigg]\bigg..
\end{align}
The combination of \eqref{14} and \eqref{16} yield that
\begin{align}\label{17}
&\mE[\Lambda_{V_{1}}+\Lambda_{V_{2}}]\no\\
&\leq C\{ I_{1}+|\Psi_{1}(0,0)|^{2}+ |\Psi_{2}(0,0)|^{2}
   +\mE[ |\Phi_{1}(0,0)|^{2}+ |\Phi_{2}(0,0)|^{2}+   \Xi_{\Gamma_{1}(\cdot, 0)}+ \Xi_{\Gamma_{2}(\cdot, 0)}]\}.
\end{align}
Next, applying It\^{o}'s formula to $\langle X_{1}(\cdot), Y_{1}(\cdot)\rangle, \langle X_{2}(\cdot), Y_{2}(\cdot)\rangle$ and then, by summing up the results and using monotonicity conditions in Assumptions 1,  we obtain

\begin{align}\label{18}
&\mE[\langle \Phi_{1}(X_{1}(T), X_{2}(T)), X_{1}(T) \rangle]+\mE[\langle \Phi_{2}(X_{1}(T), X_{2}(T)), X_{2}(T) \rangle]]\no\\
&\leq\mE[\langle \Psi_{1}(Y_{1}(0), Y_{2}(0)), Y_{1}(0) \rangle]+\mE[\langle \Psi_{2}(Y_{1}(0), Y_{2}(0)), Y_{2}(0) \rangle]\no\\
& \quad + \mE\bigg[\bigg. \int^{T}_{0}\langle \Gamma_{1}(t, \theta_{1}(t)), V_{1}(t)\rangle \dif t +\int^{T}_{0}\langle \Gamma_{2}(t, \theta_{2}(t)), V_{2}(t)\rangle \dif t \no\\
&\leq \mE[\langle \Psi_{1}(Y_{1}(0), Y_{2}(0)), Y_{1}(0) \rangle]+\mE[\langle \Psi_{2}(Y_{1}(0), Y_{2}(0)), Y_{2}(0) \rangle]\no\\
& \quad + \mE\bigg[\bigg. \int^{T}_{0}\langle \Gamma_{1}(t, 0), V_{1}(t)\rangle \dif t +\int^{T}_{0}\langle \Gamma_{2}(t,0), V_{2}(t)\rangle \dif t   \bigg]\bigg.\no\\
& \quad - L_{3}\mE\bigg[\bigg.  \int^{T}_{0}|h_{1}(t,B_{1}(t)^{\top}Y_{1}(t)+\tau\bar{B}_{1}(t)^{\top}\mE[Y_{1}(t)]+D_{1}(t)^{\top}Z_{1}(t))-h_{1}(t,0)|^{2}\dif t\no\\
&\quad \quad +\int^{T}_{0}|h_{2}(t,B_{2}(t)^{\top}Y_{2}(t)+\tau\bar{B}_{2}(t)^{\top}\mE[Y_{2}(t)]+D_{2}(t)^{\top}Z_{2}(t))-h_{2}(t,0)|^{2}\dif t    \bigg]\bigg..
\end{align}
where $\theta_{i}(t):=((\mE[X_{1}(t)])^{\top}, (\mE[Y_{1}(t)])^{\top}, (\mE[X_{2}(t)])^{\top}, (\mE[Y_{2}(t)])^{\top},   X_{i}(t)^{\top}, Y_{i}(t)^{\top}, Z_{i}(t)^{\top})^{\top}, i=1,2 . $
Applying the monotonicity conditions in  Assumption 1,   we have
\begin{align}\label{19}
&\mE[\langle \Phi_{1}(0, X_{2}(T)), X_{1}(T) \rangle]+\mE[\langle \Phi_{2}(X_{1}(T),0), X_{2}(T) \rangle]\no\\
&\leq \mE[\langle \Psi_{1}(0, Y_{2}(0)), Y_{1}(0) \rangle]+\mE[\langle \Psi_{2}(Y_{1}(0), 0), Y_{2}(0) \rangle]\no\\
& \quad + \mE\bigg[\bigg. \int^{T}_{0}\langle \Gamma_{1}(t, 0), V_{1}(t)\rangle \dif t +\int^{T}_{0}\langle \Gamma_{2}(t,0), V_{2}(t)\rangle \dif t \bigg]\bigg.\no\\
& \quad - L_{3}I_{1}.
\end{align}
By solving the above inequality and  substituting $I_{1}$ in to $\eqref{17},$    it holds that
\begin{align}\label{20}
&\mE[\Lambda_{V_{1}}+\Lambda_{V_{2}}]\no\\
&\leq C\{ |\Psi_{1}(0,0)|^{2}+ |\Psi_{2}(0,0)|^{2}
   +\mE[ |\Phi_{1}(0,0)|^{2}+ |\Phi_{2}(0,0)|^{2}+   \Xi_{\Gamma_{1}(\cdot, 0)}+ \Xi_{\Gamma_{2}(\cdot, 0)}]\}\no\\
&+\quad\mE[\langle \Psi_{1}(0, Y_{2}(0)), Y_{1}(0) \rangle]+\mE[\langle \Phi_{2}(Y_{1}(0), 0), Y_{2}(0) \rangle]]\no\\
& \quad + \mE\bigg[\bigg. \int^{T}_{0}\langle \Gamma_{1}(t, 0), V_{1}(t)\rangle \dif t +\int^{T}_{0}\langle \Gamma_{2}(t, 0), V_{2}(t)\rangle \dif t \no\\
&\quad -\mE[\langle \Phi_{1}(0, X_{2}(T)), X_{1}(T) \rangle]-\mE[\langle \Phi_{2}(X_{1}(T),0), X_{2}(T) \rangle].
\end{align}
By  a standard calculus(noting that $\epsilon,\varepsilon$ are small enough),  we derived
\begin{align}\label{20}
&\mE[\Lambda_{V_{1}}+\Lambda_{V_{2}}]\no\\
&\leq C_{1}\{ |\Psi_{1}(0,0)|^{2}+ |\Psi_{2}(0,0)|^{2}
 + \mE[ |\Phi_{1}(0,0)|^{2}+ |\Phi_{2}(0,0)|^{2}+ \Xi_{\Gamma_{1}(\cdot, 0)} +\Xi_{\Gamma_{2}(\cdot, 0)} ]\}\no\\
&\quad +\frac{1}{2}\mE[\Lambda_{V_{1}}+\Lambda_{V_{2}}].
\end{align}
We derive the desired result for \eqref{11}.

Next, we intend to prove the second result \eqref{12}.  Set
$$\grave{\Psi}_{i}(y_{1},y_{2}):=\Psi_{i}(y_{1}+\tilde{Y}_{1}(0), y_{2}+\tilde{Y}_{2}(0))-\tilde{\Psi}_{i}(\tilde{Y}_{1}(0), \tilde{Y}_{2}(0)), i=1,2,$$
$$\grave{\Phi}_{i}(x_{1},x_{2}):=\Phi_{i}(x_{1}+\tilde{X}_{1}(T), x_{2}+\tilde{X}_{2}(T))-\tilde{\Phi}_{i}(\tilde{X}_{1}(T), \tilde{X}_{2}(T)), i=1,2, $$
$$\grave{\Gamma}_{i}(t,\theta):=\Gamma_{i}(t,\theta+\tilde{\theta}_{i}(t))-\tilde{\Gamma}_{i}(t,\tilde{\theta}_{i}(t)), i=1,2, $$
for any $\theta\in \mR^{n+n+n+n+n+n+nd}.$ One can readily verify that Assumption 1  remains valid under the new set of coefficients $(\grave{\Psi}_{i},\grave{\Phi}_{i},\grave{\Gamma}_{i})$ with the same constants, $H, B_{i}, \bar{B}_{i}, D_{i},i=1,2$ in Assumption 1  and the new adjoint functions:

$$\bar{\grave{h}}_{11}(v):=\bar{h}_{11}(v+\frac{H^{\top}Y_{1}(0)-\tau H^{\top}Y_{2}(t)}{1-\tau^{2}}),$$
$$\bar{\grave{h}}_{12}(v):=\bar{h}_{12}(v+\frac{H^{\top}Y_{2}(0)-\tau H^{\top}Y_{1}(0)}{1-\tau^{2}}),$$
$$\bar{\grave{h}}_{21}(v):=\bar{h}_{21}(v+\frac{H^{\top}Y_{2}(0)-\tau H^{\top}Y_{1}(t)}{1-\tau^{2}}),$$
$$\bar{\grave{h}}_{22}(v):=\bar{h}_{22}(v+\frac{H^{\top}Y_{1}(0)-\tau H^{\top}Y_{2}(0)}{1-\tau^{2}}),$$
$$\grave{h}_{1}(t,u):=h_{1}(t,u+B_{1}(t)^{\top}Y_{1}(t)+\tau\bar{B}_{1}(t)^{\top}\mE[Y_{1}(t)]+D_{1}(t)^{\top}Z_{1}(t)),$$
$$\grave{h}_{2}(t,u):=h_{2}(t,u+B_{2}(t)^{\top}Y_{2}(t)+\tau\bar{B}_{2}(t)^{\top}\mE[Y_{2}(t)]+D_{2}(t)^{\top}Z_{2}(t)).$$
Moreover, we also confirm that the process $\hat{V}_{i}(\cdot) := V_{i}(\cdot) - \tilde{V}_{i}(\cdot) \in \mM_{\mathbb{F}}^2(\mathbb{R}^{n+n+nd})$ fulfills the FBSDE  with $(\grave{\Psi}_{i}, \grave{\Phi}_{i}, \grave{\Gamma}_{i})$. Consequently, the estimate \eqref{11} applied to $\hat{V}_{i}(\cdot)$ results in \eqref{12}.

\end{proof}

\begin{proof}[The proof of \ref{t4}]
Let $\delta_0 > 0$ be specified below, and let $\delta \in (0, \delta_0 \wedge (1 - \alpha_0)]$. Define $\alpha = \alpha_0 + \delta$ and $(\xi_{i}, \zeta_{i}, \beta_{i}(\cdot)) \in \mathbb{R}^n \times L^2_{\mathcal{F}_T}(\mathbb{R}^n) \times \sM^2_{\mathbb{F}}(\mathbb{R}^{n+n+nd})$. For any $V_{i}(\cdot) = (X_{i}(\cdot)^{\top}, Y_{i}(\cdot)^{\top}, Z_{i}(\cdot)^{\top})^{\top} \in \mM^2_{\mathbb{F}}(\mathbb{R}^{n+n+nd})$,  set
$\theta(\cdot):=(\theta_{1}(\cdot)^{\top},\theta_{2}(\cdot)^{\top})^{\top}$,
where $$\theta_{1}(\cdot):=((\mE[X_{1}(\cdot)])^{\top},(\mE[Y_{1}(\cdot)])^{\top},(\mE[X_{2}(\cdot)])^{\top},
(\mE[Y_{2}(\cdot)])^{\top},X_{1}(\cdot)^{\top},Y_{1}(\cdot)^{\top},
Z_{1}(\cdot)^{\top})^{\top},$$
$$\theta_{2}(\cdot):=((\mE[X_{1}(\cdot)])^{\top},(\mE[Y_{1}(\cdot)])^{\top},(\mE[X_{2}(\cdot)])^{\top},
(\mE[Y_{2}(\cdot)])^{\top},X_{2}(\cdot)^{\top},Y_{2}(\cdot)^{\top},
Z_{2}(\cdot)^{\top})^{\top}.$$
We consider the following MF-FBSDE:
\begin{align}\label{25}
\begin{cases}
&\dif \bar{X}_{1}(t)=[b_{1}^{\alpha}(t,\Theta_{1}(t))+\bar{\psi}_{1}(t)]\dif t
+\sum^{d}_{j=1}[\sigma_{1j}^{\alpha_{0}}(t,\Theta_{1}(t))+\bar{\gamma}_{1j}(t)]\dif W_{j}(t), t\in [0,T]\\
&\dif \bar{Y}_{1}(t)(t) =\{[f_{1}^{\alpha_{0}}(t,\Theta_{1}(t))+\bar{\phi}_{1}(t)]\dif t+\sum^{d}_{j=1} \bar{Z}_{1j}(t)\dif W_{j}(t)\\
& \bar{X}_{1}(0)=\Psi_{1}^{\alpha_{0}}(\bar{Y}_{1}(0),\bar{Y}_{2}(0))+ \bar{\xi}_{1},  \bar{Y}_{1}(T)=\Phi_{1}^{\alpha}(\bar{X}_{1}(T),\bar{X}_{2}(T) )+\bar{\zeta}_{1},
\end{cases}
\end{align}

 \begin{align}\label{26}
\begin{cases}
&\dif \bar{X}_{2}(t)=[b_{2}^{\alpha_{0}}(t,\Theta_{2}(t))+\bar{\psi}_{2}(t)]\dif t
 +\sum^{d}_{j=1}[\sigma_{2j}(t,\Theta_{2}(t))+\bar{\gamma}_{2j}(t)]\dif W_{j}(t), t\in [0,T]\\
&\dif \bar{Y}_{2 }(t) =[f_{2}^{\alpha_{0}}(t,\Theta_{2}(t))+\bar{\phi}_{2}(t)]\dif t+\sum^{d}_{j=1} \bar{Z}_{2j}(t)\dif W_{j}(t)\\
& \bar{X}_{2 }(0)=\Psi_{2}^{\alpha_{0}}(\bar{Y}_{1}(0),\bar{Y}_{2}(0))+ \bar{\xi}_{2},  \bar{Y}_{2}(T)=\Phi_{2}^{\alpha_{0}}(\bar{X}_{1}(T),\bar{X}_{2}(T) )+\bar{\zeta}_{2},
\end{cases}
\end{align}
where
\begin{align*}
&\bar{\xi}_{i}:=\xi_{i}+\delta[\Psi_{i}(Y_{1}(0),Y_{2}(0))-\Psi^{0}_{i}(Y_{1}(0),Y_{2}(0))],\\
&\bar{\zeta}_{i}:=\zeta_{i}+\delta[\Phi_{i}(X_{1}(T),X_{2}(T))-\Phi^{0}_{i}(X_{1}(T),X_{2}(T))],\\
&\bar{\beta}_{i}(\cdot):=\beta_{i}(\cdot)+\delta[\Gamma_{i}(\cdot,\theta_{i}(\cdot))-\Gamma^{0}_{i}(\cdot,\theta_{i}(\cdot))].\\
\end{align*}
One can easily check that $(\bar{\xi}_{i}, \bar{\zeta}_{i}, \bar{\beta}_{i}) \in \mR^{n}\times L^{2}_{\sF_{T}}(\mR^{n})\times \sM_{\mF}^{2}(\mR^{n+n+nd}), i=1,2.$  From the assumptions in lemma, it holds that Eqs.\eqref{25} and \eqref{26} admit a unique solution $\Theta(\cdot):=(\Theta_{1}(\cdot)^{\top},\Theta_{2}(\cdot)^{\top})^{\top},$ $$\Theta_{1}(\cdot):=((\mE[\bar{X}_{1}(\cdot)])^{\top},(\mE[\bar{Y}_{1}(\cdot)])^{\top},(\mE[\bar{X}_{2}(\cdot)])^{\top},
(\mE[\bar{Y}_{2}(\cdot)])^{\top},\bar{X}_{1}(\cdot)^{\top},\bar{Y}_{1}(\cdot)^{\top},
\bar{Z}_{1}(\cdot)^{\top})^{\top},$$
$$\Theta_{2}(\cdot):=((\mE[\bar{X}_{1}(\cdot)])^{\top},(\mE[\bar{Y}_{1}(\cdot)])^{\top},(\mE[\bar{X}_{2}(\cdot)])^{\top},
(\mE[\bar{Y}_{2}(\cdot)])^{\top},\bar{X}_{2}(\cdot)^{\top},\bar{Y}_{2}(\cdot)^{\top},
\bar{Z}_{2}(\cdot)^{\top})^{\top}.$$
Thus, given the arbitrariness of $\theta(\cdot):=(\theta_{1}(\cdot)^{\top},\theta_{2}(\cdot)^{\top})^{\top}$, we proceeded to define a mapping $G:$ $$\mM_{\mF}^{2}(\mR^{n+n+nd})\times \mM_{\mF}^{2}(\mR^{n+n+nd}) \rightarrow \mM_{\mF}^{2}(\mR^{n+n+nd})\times \mM_{\mF}^{2}(\mR^{n+n+nd}), $$
$$\theta(\cdot) \rightarrow \Theta(\cdot).  $$  If we can prove that $G$ is a   contractive mapping  when $\delta$ is
small enough, we can easily  get the result in the lemma.  For given $\theta(\cdot), \tilde{\theta}(\cdot) \in \mM_{\mF}^{2}(\mR^{n+n+nd})\times \mM_{\mF}^{2}(\mR^{n+n+nd}), $  set $\Theta(\cdot):=G(\theta(\cdot)), \tilde{\Theta}(\cdot):=G(\tilde{\theta}(\cdot)), \hat{l}:=l-\tilde{l}, l:=\theta, \Theta.  $
From estimate \eqref{12} in Lemma \ref{t1}, we then obtain
\begin{align*}
&\|\hat{\Theta}_{1}(\cdot)\|_{\mM_{\mF}^2(\mathbb{R}^{n+n+nd})}+\|\hat{\Theta}_{2}(\cdot)\|_{\mM_{\mF}^2(\mathbb{R}^{n+n+nd})} = \mathbb{E}\bigg[\bigg.\Lambda_{\hat{\Theta}_{1}}\bigg]\bigg.+\mathbb{E}\bigg[\bigg.\Lambda_{\hat{\Theta}_{2}}\bigg]\bigg.\\
&\leq \delta^2 C\bigg\{\bigg.|\check{\Psi}_{1} - \check{\Psi}_{1}^0|^2 + \mathbb{E}\bigg[\bigg.|\check{\Phi}_{1}|^2 + \bigg(\bigg.\int_0^T |\check{f}_{1}(t)|dt\bigg)\bigg.^2 \bigg]\bigg.\\
&\qquad+\bigg(\bigg.\int_0^T |\check{b}_{1}(t) - \check{b}^0_{1}(t)|dt\bigg)\bigg.^2 + \int_0^T |\check{\sigma}_{1}(t) - \check{\sigma}^0_{1}(t)|^2 dt\\
&\quad +\bigg\{\bigg.|\check{\Psi}_{2} - \check{\Psi}^0_{2}|^2 + \mathbb{E}\bigg[\bigg.|\check{\Phi}_{2}|^2 + \bigg(\bigg.\int_0^T |\check{f}_{2}(t)|dt\bigg)\bigg.^2 \bigg]\bigg.\\
&\qquad+\bigg(\bigg.\int_0^T |\check{b}_{2}(t) - \check{b}^0_{2}(t)|dt\bigg)\bigg.^2 + \int_0^T |\check{\sigma}_{2}(t) - \check{\sigma}^0_{2}(t)|^2 dt\bigg\}\bigg.,
\end{align*}

\begin{align*}
    &\check{\Psi}_{i}:= \Psi_{i}(Y_{1}(0), Y_{2}(0)) - \Psi_{i}(\tilde{Y}_{1}(0), \tilde{Y}_{2}(0)), \quad \check{\Psi}_{i}^{0} := \Psi^{0}_{i}(Y_{1}(0), Y_{2}(0)) - \Psi^{0}_{i}(\tilde{Y}_{1}(0), \tilde{Y}_{2}(0)) \\
    &\check{\Phi}_{i}:= \Phi_{i}(X_{1}(T), X_{2}(T)) - \Phi_{i}(\tilde{X}_{1}(T), \tilde{X}_{2}(T)) \\
    &\check{l}_{i}(t):= l_{i}(t, \theta_{i}(t)) - l_{i}(t, \tilde{\theta}_{i}(t)) \quad \text{with} \quad l_{i} = f_{i}, b_{i}, b^0_{i}, \sigma_{i}, \sigma^0_{i}, i=1,2.
\end{align*}
Basing on the Lipschitz condition, it holds that

$$\| \hat{\Theta}(\cdot) \|_{\mM_{\mathbb{F}}^2(\mathbb{R}^{n+n+nd})\times \mM_{\mathbb{F}}^2(\mathbb{R}^{n+n+nd})}^2 \leq \delta^2 C \| \hat{\theta}(\cdot) \|_{\mM_{\mathbb{F}}^2(\mathbb{R}^{n+n+nd})\times \mM_{\mathbb{F}}^2(\mathbb{R}^{n+n+nd})}^2.$$
Choosing $\delta_{0}:=\frac{1}{2\sqrt{C}},$ then for any $\delta\in (0,\delta_{0}\wedge (1-\alpha_{0})],$ the above inequality shows that $G$ is contractive.
It is evident that the unique fixed point corresponds precisely to the unique solution of System $(\pi_{1})$  when $\alpha=\alpha_{0}+\delta$ and $\xi_{i},\zeta_{i},\beta_{i}(\cdot), i=1,2.$ The proof is complete.

\end{proof}
By applying the above two lemmas, we give the proof of Theorem \ref{tt}].
\begin{proof}[The proof of Theorem \ref{tt}]
As noted earlier, System $(\pi_{1})$  is uniquely solvable when $\alpha_{0}=0.$ By iteratively applying Lemma \ref{t4}, the unique solvability is then extended from $\alpha = 0$ to $\alpha > 0.$ Given that the step size $\delta_{0}>0$ is fixed, only a finite number of such extensions are needed to establish the unique solvability of System $(\pi)$ for $\alpha=1$, thereby completing the proof.

\end{proof}

\subsection{Convexity and Uniform Convexity}
We list some auxiliary results in this section.
The following results can be found in \cite{Peypouquet2015,LiuNiuWangYu2026}.
\bl\label{A1}
 Under Assumption 3, both $\nabla f_{1i}(\cdot): \mathbb{R}^{m} \rightarrow \mathbb{R}^{m}, i=1,2$ and $\nabla f_{4i}(\cdot): \mathbb{R}^{k} \rightarrow \mathbb{R}^{k},i=1,2$ are bijective. Let $(\nabla f_{1i})^{-1}(\cdot), i=1,2$ and $(\nabla f_{4i})^{-1}( \cdot), i=1,2$ denote the four inverse mappings, respectively. Then they are Lipschitz continuous with Lipschitz constant $1 / \delta>0$. Moreover, $(\nabla f_{4i})^{-1}$ is $ \sB(\mathbb{R}^{k})$-measurable.
\el

\bd \label{d1}
Suppose $D \subset \mathbb{R}^n$ is a nonempty and convex set. A function $g : D \rightarrow \mathbb{R}$ is convex if
$$
g(\lambda x + (1 - \lambda)y) \leq \lambda g(x) + (1 - \lambda) g(y)
$$
for any $\lambda \in (0, 1)$ and any $x, y \in D$. If the inequality holds strictly whenever $x \neq y$, then $f(\cdot)$ is called strictly convex. Furthermore, $g(\cdot)$ is uniformly convex (also known as strongly convex) with parameter $\delta > 0$ if
$$
g(\lambda x + (1 - \lambda)y) + \frac{\delta}{2} \lambda (1 - \lambda) |x - y|^2 \leq \lambda g(x) + (1 - \lambda) f(y)
$$
for any $\lambda \in (0, 1)$ and any $x, y \in D$.
\ed

 \bl \label{100}
 Let $D \subset \mathbb{R}^n$ be a nonempty, open, and convex set.  Suppose $g : D \rightarrow \mathbb{R}$ is differentiable.  Then, the following statements are equivalent.

\begin{enumerate}
    \item [(1)]$g(\cdot)$ is convex (respectively, uniformly convex with $\delta > 0$).
    \item[(2)] $g(x) - g(y) - \langle \nabla g(y), x - y \rangle \geq 0$ (respectively, $\geq (\delta/2) |x - y|^2$) for any $x, y \in D$.
    \item[(3)] $\langle \nabla g(x) - \nabla g(y), x - y \rangle \geq 0$ (respectively, $\geq \delta |x - y|^2$) for any $x, y \in D$.
\end{enumerate}

\el

\subsection{Projection Onto a Closed Convex Set}
In this section, we list some basic properties of a projection onto a closed convex subset of $\mathbb{R}^n$. More details can be found in \cite{Brezis2011}. With a slight abuse of notation, the inner product and the induced norm of $\mathbb{R}^n$ (not necessarily the Euclidean inner product and the Euclidean norm) are denoted by $\langle \cdot, \cdot \rangle$ and $|\cdot|$, respectively.

\bl\label{c}
 Let $K \subset \mathbb{R}^n$ be a nonempty closed convex set. Then, for each $x \in \mathbb{R}^n$, there exists a unique element $\Pi(x) \in K$ such that
\[
|x - \Pi(x)| = \min_{y \in K} |x - y|.
\]
\el
\noindent
The element $\Pi(x)$ is called the projection of $x$ onto $K$. Moreover, $\Pi(x) \in K$ is characterized by the property that
\begin{equation}\label{a}
\langle x - \Pi(x), y - \Pi(x) \rangle \leq 0
\end{equation}
for any $y \in K$. Furthermore,
\begin{equation}\label{b}
\begin{cases}
|\Pi(x) - \Pi(\bar{x})|^2 \leq \langle \Pi(x) - \Pi(\bar{x}), x - \bar{x} \rangle \\
|\Pi(x) - \Pi(\bar{x})| \leq |x - \bar{x}|
\end{cases}
\end{equation}
for any $x, \bar{x} \in \mathbb{R}^n$.

\section*{Declaration of competing interest}
The authors declare that they have no known competing financial interests or personal relationships that could have appeared
to influence the work reported in this paper.

\section*{Funding}
This research is supported by the National Natural Science Foundation of China (Grant no.  11626236), the Fundamental Research Funds for the Central Universities of South-Central Minzu University (Grant nos. CZY15017).

\end{document}